\documentclass[a4paper]{article}
\usepackage[english]{babel}
\usepackage{pgfplots}\pgfplotsset{compat=1.18}
\usepackage{float}
\usepackage[T1]{fontenc}
\usepackage[utf8]{inputenc}
\usepackage[affil-it]{authblk}
\usepackage{algpseudocode}
\usepackage{paralist}
\usepackage{caption}
\usepackage{appendix}
\usepackage{verbatim}
\usepackage{wrapfig}
\usepackage{amsmath}
\usepackage[arrow, matrix, curve]{xy}
\usepackage{mathtools}

\usepackage[x11names]{xcolor}
\usepackage{amssymb}
\usepackage{wasysym}
\usepackage{listings}
\usepackage[gen]{eurosym}
\usepackage{bbm}
\usepackage{stmaryrd}
\usepackage{graphicx}
\usepackage{mathrsfs}
\usepackage{pst-func}
\newcommand{\dd}{\mathrm{d}}

\DeclareMathOperator{\Sym}{\mathcal{S}_{2d-2}}
\DeclareMathOperator{\Symm}{\mathcal{S}^{\mathrm{Sym}}_{2d-2}}
\DeclareMathOperator{\Symtwo}{\R^{2\times 2}_{\mathrm{Sym}}}

\DeclareMathOperator{\tr}{\mathrm{tr}}

\DeclareMathOperator{\id}{\mathrm{id}}
\DeclareMathOperator{\R}{\mathbb{R}}
\DeclareMathOperator{\C}{\mathbb{C}}
\DeclareMathOperator{\G}{\mathbb{G}}
\DeclareMathOperator{\N}{\mathbb{N}}
\DeclareMathOperator{\diag}{\mathrm{diag}}

\usepackage{enumitem}
\usepackage{tabularx}
\usepackage{amsthm}

\usepackage{tikz}
\usetikzlibrary{decorations, decorations.text}
\usetikzlibrary{plotmarks}
\newtheorem{Theorem}{Theorem}[section]
\newtheorem{corollary}[Theorem]{Corollary}
\newtheorem{lemma}[Theorem]{Lemma}
\newtheorem{Bem}[Theorem]{Remark}
\newtheorem{Proposition}[Theorem]{Proposition}
\newtheorem{Definition}[Theorem]{Definition}
\newtheorem{Scheme}[Theorem]{Scheme}

\usepackage{fancyhdr}
\usepackage{datetime}
\renewcommand{\sectionmark}[1]{\markboth{~#1}{}}

\fancypagestyle{main}{
\fancyhf{}
\fancyhead[LE]{\footnotesize\nouppercase{\sectionmark}}
\fancyhead[RO]{\footnotesize \nouppercase{\sectionmark}}
\fancyfoot[LE,RO]{\thepage}

}
\renewcommand\subsectionmark[1]{}

\usepackage{csquotes}
\usepackage[backend = biber, sorting=nyt, style =numeric-comp, giveninits=true,
  maxbibnames=99]{biblatex}

\usepackage[breaklinks]{hyperref}

\begin{document}

\title{\bfseries Sum-of-squares rank-one convexity: a semidefinite approach to Morrey's problem}

\author[1]{{Florian Gammer}}

\affil[1]{\footnotesize Faculty of Mathematics, University of Regensburg, Universitätsstraße 31
93053 Regensburg, Germany}

\date{}

\maketitle

\begin{abstract}
Morrey's problem asks whether rank-one convex functions are automatically quasiconvex. 
Here, we introduce the class of sum-of-squares rank-one convex polynomials and show that every rank-one convex function on $\R^{m\times n}$ can be locally uniformly approximated by such polynomials if and only if $\min\{m,n\}\le 2$. These polynomials provide a tractable and constructible subclass and the density reduces Morrey's problem to the study of this subclass when $\min\{m,n\}\le 2$.

In dimension $2\times2$, we characterize these polynomials by decomposing their Hessians into a sum of squares and a polynomial multiple of the determinant. This yields an exact semidefinite parameterization for each fixed degree. For every prescribed smooth periodic test function $u$, a rank-one convex integrand violates the quasiconvexity inequality at the origin with witness $u$ if and only if the resulting semidefinite program detects a violation by the same $u$ at some finite degree. Whenever the normalized program is feasible, it satisfies strong duality and its dual optimum is attained.

We also prove density results for symmetric rank-one convex integrands on $\Symtwo$ and for isotropic rank-one convex integrands on $\R^{2\times2}$. In the symmetric setting, the Hessian representation admits a sum-of-squares multiplier for the negative determinant. In the isotropic setting, rotational invariance yields a block-diagonal sum-of-squares representation.
\end{abstract}

%90C23: Polynomial optimization
%90C22: Semidefinite programming
\smallskip


\section{Introduction}

The problem of minimizing an energy
\begin{equation}\label{eq: general energy}
    \mathbb E(u) = \int_\Omega f(Du) \dd x,
\end{equation}
where $\Omega \subset \mathbb R^n$ is a bounded Lipschitz domain, $f \colon \mathbb R^{m \times n} \to \mathbb R$ is a continuous energy density, and $u \in W_g^{1,p}(\Omega, \mathbb R^m)$ is a Sobolev map with $g \in W^{1,p}(\Omega, \mathbb R^m)$, lies at the heart of the calculus of variations. A minimizer is guaranteed to exist if the energy functional is coercive and lower semicontinuous. 

Under suitable growth conditions, the lower semicontinuity of $\mathbb E$ is equivalent to the quasiconvexity of $f$ \cite{morrey1952quasi}, which means that for all $A \in \mathbb R^{m \times n}$ and all $v \in C^\infty(\mathbb T^n, \mathbb R^m)$,
$$f(A) \le \int_{\mathbb T^n} f(A + D v) \dd x.$$
Quasiconvexity turns out to be highly intractable due to its inherent nonlocality (cf.~\cite{kristensen1999non}). Any quasiconvex function also satisfies the local property of rank-one convexity. A continuous function is called rank-one convex if, for all $A, B \in \mathbb R^{m \times n}$ with $\mathrm{rk}(A-B) = 1$ and all $t \in (0,1)$,
$$f(A + t(B-A)) \le (1-t)f(A) + t f(B).$$
These two notions of convexity are very closely related, and the question of whether they coincide is known as Morrey's problem.

Although equivalence has been established for certain subclasses of the rank-one convex functions (see, e.g., \cite{guerra2026quasiconvexityburkholderfunctionsymmetric,KristensenGuerra2022,chenchiah2008relaxation,grabovsky2019explicit,grabovsky2019rank,kirchheimautomatic,conti2003polyconvexity,harris2018two,muller1999rank,martin2017rank,faraco2003quasiconvex,bandeira2009finding}),
the converse implication fails in several matrix dimensions. The first such counterexample was obtained by \v Sverák \cite{vsverak1992rank} in $\R^{3\times 2}$, and Grabovsky
\cite{grabovsky2018microstructure} later obtained another one in $\mathbb R^{8\times2}$. 
More recently, Cassese
\cite{cassese2026solutionmorreysproblemmathbbr2times}
constructed a homogeneous quartic counterexample in
$\mathbb R^{2\times3}$ using martingale theory, leaving $\mathbb R^{2\times2}$
as the only unresolved dimension
for the general Morrey problem (for more details and further counterexamples in higher dimensions,  cf.~\cite{ConDolDon2026,gammer2026symmetricmorreyproblemdimensions,ball1987does,agazzi2026morreysproblemmathbbr2times} and the references therein).

Since no analytical counterexample has yet been constructed in $\R^{2\times 2}$, researchers have also turned to numerical schemes (see, e.g., \cite{GUTIERREZ20073177,guerra2021numerical,voss2022numerical}). Such approaches face two complementary challenges: constructing sufficiently general rank-one convex integrands while identifying periodic test functions capable of detecting a violation of quasiconvexity. Consequently, numerical investigations necessarily rely on restricted classes of integrands and test functions, and the absence of a detected counterexample may reflect either the validity of Morrey's conjecture or limitations of the considered rank-one convex functions $f$ or test functions $u$.

 In this context, the density results of Heinz \cite{Heinz2008Preservation} are particularly noteworthy. They establish that every rank-one convex (resp.~quasiconvex) function can be locally uniformly approximated by rank-one convex (resp.~quasiconvex) polynomials. This allows us to restrict the search to polynomial integrands, providing a more tractable class of functions to investigate. A natural test for Morrey's conjecture would then mirror \v Sverák's example by fixing a test function $u$ and varying over rank-one convex polynomials. 
A key bottleneck for employing this approach is identifying rank-one convex polynomials. 
Indeed, even verifying whether a given quartic polynomial on $\R^{2\times n}$ is rank-one convex is strongly NP-hard. This follows from the strong NP-hardness of deciding convexity of quartic polynomials established in~\cite{ahmadi2013np}, since ordinary convexity can be embedded into rank-one convexity by considering integrands of the form $p(A)=f(A^Te_1)$.

In this paper, we use techniques from optimization and real algebraic geometry to identify a constructible dense subclass of the rank-one convex polynomials. Verifying their defining property is a finite-dimensional semidefinite feasibility test (cf.~\cite{vandenbergheisemidefinite}). Moreover, this approach yields an explicit finite-dimensional parametrization of all such polynomials of a prescribed maximal degree, enabling a new numerical approach to Morrey's problem that is sufficient, in the sense that any counterexample to the general problem leads to a counterexample in the constructible subclass. For the appropriate choice of test function and a sufficiently large degree bound, such a counterexample will be contained in the admissible SDP class. 

Let us provide some background on the algebraic tools we use. We can reformulate the rank-one convexity condition for a polynomial function $p\in \mathcal{P}(\R^{m\times n})$ in terms of its second derivative. Then $p$ is rank-one convex if and only if, for all $A\in \R^{m\times n}$, $a\in \R^m$, and $b\in \R^n$, there holds
\begin{equation}
    D^2p(A)[a\otimes b,a\otimes b]\ge 0 \label{eq: non-negativity}.
\end{equation}
Crucially, $P=D^2(\cdot)[\cdot,\cdot]$ is a polynomial, and \eqref{eq: non-negativity} requires it to be nonnegative. Artin's Theorem \cite{artin1927zerlegung} establishes that this is equivalent to the existence of real polynomials $r_1,\dots,r_N,q\in \mathcal{P}(\R^{m\times n}\times \R^m\times \R^n)$, $q\ne 0$, such that 
$$q^2P=\sum_{i=1}^Nr_i^2.$$
Nonnegative polynomials that admit a decomposition with $q=1$ are called \emph{sum-of-squares} (SOS) polynomials. Crucially, $q$ cannot be set to $1$ in general (\cite{hilbert1888darstellung,motzkin1967arithmetic,blekherman2006there}). The SOS polynomials therefore form a strict subclass (see \cite[Thm.~9.1]{berg1976positive}, along with \cite{lasserre2007sum,lassere2007sos}) of the nonnegative polynomials, and as it turns out, they are dense therein in the locally uniform topology. 
These methods have become an important tool in optimization, and have recently gained interest in the analysis and PDE communities. For instance, it was used to compute quasiconvex envelopes and study polyconvex functions (\cite{fantuzzi2026polyconvexitymomentssumssquares,henrion2026semidefiniterelaxationsnonlinearelasticity}). The definition of SOS polyconvexity leads us to the main notion of this paper.
\begin{Definition}
    A polynomial $p\in \mathcal{P}(\R^{m\times n})$ is \emph{sum-of-squares rank-one convex} ($\mathrm{SOSRC}$) if there are polynomials $r_1,...,r_N\in \mathcal{P}(\R^{m\times n}\times \R^m\times\R^n)$, $N\in \N$, such that for all $A\in \R^{m\times n}$, $a\in \R^m$, $b\in \R^n$ there holds
    \begin{equation}\label{eq: SOSRC condition}
         D^2p(A)[a\otimes b,a\otimes b]=\sum_{\nu=1}^Nr_{\nu}(A,a,b)^2.
    \end{equation}
\end{Definition}
Note that finding an SOS representation pointwise in $A$ is always possible (at least when $\min\{m,n\}=2$) by a variant of Hilbert's Theorem (cf.~\cite{bandeira2016role}). Our definition of SOSRC polynomials, by contrast, requires a \emph{simultaneous} SOS representation in all variables $A,a,b$. This is a much stronger condition,
for instance, take a convex bivariate sextic $h$ that is not SOS-convex (see \cite{ahmadi2013np}), and set $p(A)=h(A_{11},A_{12})$. Then $p$ is rank-one convex but not SOSRC\footnote{Otherwise, $p$, restricted to $a=e_1$ and $A_{21}=A_{22}=0$, would give an SOS representation of $b^TD^2h(x)b$}.
Moreover, since verifying the SOS condition is a semidefinite feasibility problem rather than a general polynomial nonnegativity problem. Finally, this class of functions admits an explicit finite-dimensional parametrization, as we will discuss below.

\subsubsection*{Main results I: Density}
Although rank-one convex polynomials are dense among rank-one convex integrands, it is not immediate that one can impose an SOS condition for the Hessian on rank-one matrices without losing density. Our first result establishes precisely this stronger approximation property when $\min\{m,n\}=2$. In the case $\min\{m,n\}<2$ rank-one convexity is just convexity, and the density of SOS-convex polynomials is well-established (cf., e.g, \cite{lassere2007sos}).
\begin{Theorem}\label{Thm: Density of SOSRC}
Let $\min\{m,n\}=2$.  
    Let $f\in C^0(\R^{m\times n})$ be a rank-one convex function. Then for any compact $K\subset  \R^{m\times n} $ and any $\varepsilon>0$, there is an $\mathrm{SOSRC}$ polynomial $p\in   \mathcal{P}(\R^{m\times n})$ such that 
    $$\sup_{A\in K}\vert f(A)-p(A)\vert <\varepsilon.$$
\end{Theorem}
We also show that the density result is false whenever $m,n\ge 3$. 
The immediate consequence of the density theorem is that it restricts Morrey's problem to the class of SOSRC polynomials.
\begin{corollary}\label{Cor: SOSRC Morrey} In dimension $2\times2$, the following are equivalent:
   \begin{itemize}
       \item[(i)] Every $\mathrm{SOSRC}$ polynomial $p$ is quasiconvex.
 \item[(ii)]  Every continuous rank-one convex integrand $f$ is quasiconvex.
   \end{itemize}
\end{corollary}

\subsubsection*{Main results II: Representation}
The SOSRC condition specifies the Hessian only on rank-one directions. In dimension $2\times2$, every quadratic form vanishing on rank-one matrices is a multiple of the determinant. We can therefore extend the SOS representation to arbitrary matrix directions by adding a polynomial multiple $\Lambda(X)\det H$. This multiplier must be chosen to be compatible with the Hessian structure. For notational convenience, we identify $\R^{2\times 2}\cong \R^4$ via $X\mapsto (X_{11},X_{12},X_{21},X_{22})$.
\begin{Theorem}\label{Thm: Representation}
    Let $p\in \mathcal{P}(\R^4)$ be an $\mathrm{SOSRC}$ polynomial of degree at most $2d$. Then there is a decomposition \begin{equation} \label{eq: split of SOSRC polynomials}p=p_{\mathrm{aff}}+p_{\mathrm{sos}}+p_{\mathrm{det}},
    \end{equation}
where $p_{\mathrm{aff}}$ is affine, $p_{\mathrm{sos}}$ is a sum-of-squares polynomial of degree at most $2d$, and $p_{\mathrm{det}}$ is the product of the determinant and a polynomial of degree at most $2d-2$.\\
More precisely, we can construct $p$ from the SOSRC representation of $D^2p$ on rank-one matrices. There exist an integer $N\in \N$, a polynomial $\Lambda\in \mathcal{P}(\R^4)$ of degree $\deg \Lambda\le 2d-2$, and polynomials
    $v_{\nu}\in \mathcal{P}(\R^4)^4$ of degree $\deg v_{\nu,i}\le d-1 $ for $\nu=1,...,N$, $i\in \{1,...,4\}$ such that for $v_\nu(X)=\sum_{\vert\alpha\vert\le d-1}v_{\nu,\alpha}X^\alpha$, and $\Lambda(X)=\sum_{\vert\beta\vert\le 2d-2}\lambda_\beta X^\beta$, the matrix
$$Q_\gamma=\sum_{\nu=1}^N\sum_{\alpha+\beta=\gamma}
 v_{\nu,\alpha}\otimes v_{\nu,\beta}+\lambda_\gamma \mathbb{D},$$
where $\mathbb{D}$ represents the quadratic form associated with the determinant, satisfies the compatibility conditions for all $i,j,k\in \{1,...,4\}$, $\vert \gamma \vert \le 2d-3$
\begin{equation}\label{eq: coefficient compatibility}
(\gamma_k+1)(Q_{\gamma+e_k})_{ij}  =(\gamma_i+1)(Q_{\gamma+e_i})_{kj}.
\end{equation}
Let us write
$v^{\mathrm{hom},j}_{\nu}(X)=\sum_{\vert \alpha\vert = j}v_{\nu,\alpha}X^\alpha,$ (resp.~$\Lambda^{\mathrm{hom}}_j$), for the $j$-homogeneous part of $v_\nu$ (resp.~$\Lambda$). Then there are $c\in \R,B\in \R^4$ such that 
\begin{align}\label{eq: integrated p}
    \begin{aligned}
        p_{\mathrm{aff}}(X)&=c+B^TX, \qquad p_{\mathrm{det}}(X)=\biggl[\sum_{k=0}^{2d-2}
 \frac{\Lambda^{\mathrm{hom}}_k(X)}{(k+1)(k+2)}\biggr]\det X,\\
 p_{\mathrm{sos}}(X)&= \sum_{\nu=1}^N\sum_{i,j=0}^{d-1}
 X^T \frac{v^{\mathrm{hom},i}_{\nu}(X)\otimes v^{\mathrm{hom},j}_{\nu}(X)}{(i+j+1)(i+j+2)}X.
    \end{aligned}
\end{align}
\end{Theorem}
An analogous representation (with a linear combination of $2\times 2$ minors in the second part) holds in general for $\min\{m,n\}=2$. 
The representation theorem also yields an exact semidefinite parametrization of all SOSRC polynomials of prescribed maximal degree. This parametrization is the basis for the semidefinite formulation developed below.
\begin{corollary}\label{Cor: Generation of SOSRC functions}
Let $d\in \N$, let  $m_{d-1}(X)=(X^\alpha)_{\vert \alpha\vert\le d-1}$ be the monomials up to degree $d-1$ in $4$ variables. Set $z_d(X,H)=m_{d-1}(X)\otimes H$, $N_d=4\binom{d+3}{4}$. Let $\G\in \R_{\mathrm{Sym}}^{N_d\times N_d}$ be positive semidefinite and let $\Lambda\in \mathcal{P}_{2d-2}(\R^4)$, $\lambda$ its coefficient vector. Let $C\in \R^{\binom{2d+4}{4}-5}$. For coefficients $C$, we set $p_C:=\sum_{\vert \alpha\vert =2}^{2d}C_\alpha X^\alpha$, and define a feasibility operator (coefficient matching)
$$\mathscr{F}_d(C,\G,\lambda)[X,H]:=D^2p_C(X)[H,H]-z_d(X,H)^T\G z_d(X,H)-\Lambda(X)\det(H).$$
Then, a polynomial $p$ on $\R^{2\times 2}$ of degree at most $2d$ is $\mathrm{SOSRC}$ if and only if it can be represented as $p_C+p_{\mathrm{aff}}$, where $p_{\mathrm{aff}}$ is affine and there exist  $C,\G,\Lambda$  such that 
$$\mathscr{F}_d(C,\G,\lambda)\equiv 0, \quad \text{and} \quad \G\succeq 0. $$
\end{corollary}

\subsubsection*{Main results III: SDP formulation}
Combining the density result with the exact parametrization of SOSRC polynomials, we derive a semidefinite programming criterion for quasiconvexity violations. For any fixed $u\in C^\infty(\mathbb T^2,\mathbb R^2)$ and any prescribed polynomial degree $2d$, we obtain a finite-dimensional semidefinite program.
The constraint $b_u^T\lambda=-1$ normalizes the determinant contribution, which must be negative for a quasiconvexity violation. The objective then measures the SOS contribution, so a violation occurs precisely when $\rho_d(u)<1$. Finally, we may impose $\lambda_0=0$, since determinants are null Lagrangians. 
\begin{Theorem}\label{Thm: Morrey Test}
    For any $u\in C^\infty(\mathbb{T}^2,\R^2)$, and $d\in \N$, we define a matrix $\mathbb{M}_u\in \R^{N_d\times N_d}$ and a vector\footnote{We omit the constant coefficient $\lambda_0$.} $b_u\in \R^{\binom{2d+2}{4}-1}$, $\{\alpha\in\mathbb N_0^4\mid 1\le|\alpha|\le2d-2\}$ via
    \begin{align*}
        \mathbb{M}_u&=\int_{\mathbb{T}^2}\int_0^1(1-t)z_d(tDu,Du)\otimes z_d(tDu,Du)\dd t\dd x,\\
        (b_u)_{\alpha}&=\frac{1}{(\vert \alpha\vert +2)(\vert \alpha \vert +1)}\int_{\mathbb{T}^2}(Du)^\alpha\det(Du)\dd x,
    \end{align*}
    where $N_d,z_d$ are as in Cor.~\ref{Cor: Generation of SOSRC functions}. 
   Let $\mathscr{C}_d$ be the linear compatibility operator obtained by expressing the Hessian compatibility conditions in coefficient form.
We define a semidefinite program via
    \begin{align}\label{eq: semidefinite program}\tag{SDP}
        \begin{aligned}
           &\rho_d(u)=\inf_{\G,\lambda}\, \langle \mathbb{M}_u,\G\rangle \\
            &\mathrm{subject\ to} \ \mathscr{C}_d(\G,\lambda)=0, \ b_u^T\lambda= -1, \ \mathrm{and}\ \G\succeq 0.
        \end{aligned}
    \end{align}
    Then there is a rank-one convex function on $\R^{2\times 2}$ failing to be quasiconvex at $0$ with witness $u$ if and only if there is a $d\in \N$ such that $\rho_d(u)<1$.
\end{Theorem}
Thm.~\ref{Thm: Density of SOSRC} implies that Morrey's problem can be restricted to SOSRC polynomials. Thus, whenever a continuous rank-one convex integrand violates the corresponding quasiconvexity inequality, there exists a finite polynomial degree at which the associated SDP detects this violation.

The primal problem \eqref{eq: semidefinite program} is strictly feasible at fixed degree $2d$ if and only if $b_u\ne 0$. Test functions for which $b_u=0$ cannot witness a violation of quasiconvexity in the first place. In particular, strong duality holds, and the optimal value in the dual problem is attained. Moreover, the dual program contains a feasible interior point under the additional assumption that
\begin{equation}\label{eq: Strict duality condition}\tag{IP}
    \mathrm{there}\ \mathrm{are}\ \mathrm{no}\ \G\succeq 0, \G\ne 0,\ \lambda \ \mathrm{such}\ \mathrm{that}\ \mathscr{C}_d(\G,\lambda)=0, \ b_u^T\lambda= 0, \ \langle \mathbb{M}_u,\G\rangle=0.
\end{equation}
Although verifying this condition analytically becomes intractable at higher degrees, it can still be checked numerically. When \eqref{eq: Strict duality condition} holds, primal–dual interior-point methods are can be applied (at least for low degrees). 

We implemented this semidefinite program using a primal-dual interior-point method. However, we only conduct computations for low degrees ($2d=12$ at most), since the number of operations required is of order $\mathcal{O}(d^{12})$. Consequently, we could not obtain a counterexample for Morrey's conjecture. Even though the general semidefinite formulation of Morrey's problem is computationally demanding, it provides a novel framework for numerical investigations.

\subsubsection*{Main results IV: Extension to isotropic and symmetric functions}

In this final section, we derive analogous density results for symmetric rank-one convex functions on $\Symtwo$ and isotropic rank-one convex functions on $\R^{2\times 2}$. Both settings yield additional structure in the SOSRC representation, albeit in different ways. In the symmetric case, apart from reducing the number of variables from $4$ to $3$, the determinant multiplier can be chosen to be a sum-of-squares polynomial. In the isotropic case, the matrix $\G$ can be chosen to admit a certain block-diagonal structure.

For symmetric SOSRC polynomials, Finsler's lemma is no longer applicable. Instead, we employ the Yakubovich S-lemma, exploiting the fact that the symmetric rank-one directions are precisely the matrices $H$ satisfying $\det H\le 0$. This also yields the SOS representation of the determinant multiplier, since the polynomial identity, initially established on $\{\det H\le 0\}$, extends to all of $\Symtwo$ because $\{\det H<0\}$ is a nonempty open set.

For isotropic SOSRC polynomials, density within the isotropic rank-one convex functions follows by averaging the SOSRC approximating sequence with respect to the Haar measure on the group. Moreover, the block-diagonal structure of the matrix becomes apparent after complexifying the feature variables and grouping them according to their rotation labels. Since the group action is diagonal in these variables, Haar averaging eliminates interactions between features with distinct labels.

\subsubsection*{Outline}
We prove our density result in Section \ref{Sec: Density} and establish the representation and generation formulas in Section \ref{Sec: Generation}. 
In Section \ref{Sec: Morrey SDP}, we obtain the semidefinite programming formulation of Morrey's problem. Moreover, we prove strong duality, derive the criterion for dual strict feasibility, lay out our numerical scheme, and present our numerical results. 
In Section \ref{Sec: Corollaries}, we finally derive the reduced versions of the previous results for isotropic and symmetric rank-one functions.

\subsubsection*{Notation}
Throughout, we denote matrices by $A,B,H$, use $\R_\mathrm{Sym}^{N\times N}$ for the space of symmetric
$N\times N$ matrices, and write $\mathcal{P}_{d}$ for the space of real polynomials of degree up to $d$. 
We write $\mathbb I$ for the identity matrix, indicating dimensions with an index.
For integrands $f$, $D^2f(A)[B,B]$ denotes the second
directional derivative at $A$ in direction $B$.
We use $\mathbb T^n=\mathbb R^n/\mathbb Z^n$, which we identify with $(0,1)^n$, equipped with the
Lebesgue measure. We write $C^\infty(\mathbb T^n,\R^m)$ for the space of $\mathbb Z^n$-periodic smooth functions $v\colon \R^n\to \R^m$.
We denote the quadratic form representation of the determinant by
$$\mathbb{D}=\frac{1}{2}\begin{pmatrix}
    0 & 0 & 0 & 1\\
    0 & 0 & -1 & 0\\
    0& -1 & 0 & 0\\
    1& 0 & 0& 0
\end{pmatrix}.$$ 
We let
$$\Sym:=\bigl\{P\in \bigl(\mathcal{P}_{2d-2}(\R^4)\bigr)^{4\times 4}\mid P^T=P\bigr\}.$$
We use $\Lambda$ and $p_C$ to denote the corresponding polynomial functions, and $\lambda$ and $C$ for their respective coefficient vectors. For operators on polynomial spaces, we will identify the polynomial arguments with their coefficient vectors in the monomial basis.
We denote high-dimensional matrices by $\G,\mathbb{M}$. 
We denote the operator norm of a matrix by $\vert \mathbb{F}\vert_{\mathrm{op}}$.

\section{Density of SOSRC polynomials}
\label{Sec: Density}

We make great use of the following three theorems.

\begin{Proposition}[{\cite[Cor.~3.3]{lassere2007sos}}]\label{Prop: Density of SOS polynomials}
    Let $p\in\mathcal{P}(\R^{n})$ be a polynomial that is
nonnegative on the cube $Q_R:=[-R,R]^n$, $R>0$. 
For $r\in\mathbb{N}$, define the perturbation
\begin{equation}\label{eq: Theta r R}
    \Theta_{r,R}(X)
:=
1+\sum_{j=1}^n\biggl(\frac{X_j}{R}\biggr)^{2r}
\in\mathcal{P}(\R^{n}).
\end{equation}
Let $\varepsilon>0$ be fixed. Then there exists
$r_0=r_0(p,\varepsilon,R)\in\mathbb{N}$ such that, for every
$r\ge r_0$, the polynomial
$$p_{\varepsilon,r,R}
:=
p+\varepsilon\Theta_{r,R}$$
is a sum of squares.
\end{Proposition}
Lasserre and Netzer establish convergence of the coefficients in the $l_1$-norm. This implies local uniform convergence on $Q_R$, since for any compact set $K\subset Q_R$, we have
$$\sup_{x\in K}\Vert p(x)-p_{\varepsilon,r,R}(x)\Vert\le \varepsilon (1+ n).$$
\begin{lemma}\label{Lem: compact open topology}
Let $f\in C^0(\R^n)$ be given. Suppose that for all compact $K\subset \R^n$ we have
$$\inf_{p\in\mathcal{P}_{\mathrm{sosrc}}}\Vert f-p\Vert_{C^0(K)}=0.$$
Then there exists a sequence $(p_j)_{j\in\mathbb{N}}\subset\mathcal{P}_{\mathrm{sosrc}}(\R^n)$
such that
$$p_j\longrightarrow f
\quad \text{locally uniformly in }\R^n.$$
\end{lemma}
\begin{proof}
Consider the compact exhaustion
$K_j:=\overline{B_j(0)}\subset \R^n.$
By assumption, for every $j\in\mathbb{N}$, there exists
$p_j\in\mathcal{P}_{\mathrm{sosrc}}(\R^n)$ such that
$$\Vert p_j-f\Vert_{C^0(K_j)}<\frac{1}{j}.$$
Since every compact subset $K\subset \R^n$ is contained in some
$K_{j_0}$, we obtain, for all $j\ge j_0$,
$$\Vert p_j-f\Vert_{C^0(K)}
\le \Vert p_j-f\Vert_{C^0(K_j)}
<\frac{1}{j}.$$
Thus, $p_j\to f$ locally uniformly in $\R^n$.
\end{proof}
Note that we cannot use Lasserre's earlier result
(\cite[Thm.~4.1]{lasserre2007sum}), which gives an analogous construction for nonnegative polynomials on $\R^n$ that yields locally uniform convergence, because $D^2p(A)[H,H]$ need not be nonnegative for $\operatorname{rk}H\ge 2$, and we can compensate for this lack of nonnegativity (with a polynomial) only on a compact set. Nevertheless, Lemma~\ref{Lem: compact open topology} shows that taking a diagonal sequence of the Lasserre--Netzer approximations still yields a locally uniform approximation by SOS polynomials on $\R^n$.

\begin{lemma}[Finsler's lemma\footnote{Since we will apply this lemma pointwise, the reader might be tempted to suggest that we instead use the matrix polynomial version of Finsler's lemma \cite{CIMPRIC2015239}. However, this version is inapplicable in our context because the required matrix $G$ is indefinite, meaning that condition A is not satisfied.} \cite{finsler1936vorkommen}]\label{lem: Finsler's lemma}
    Let $Q,L\in\mathbb R_{\mathrm{Sym}}^{N\times N}$, where $L$ is indefinite. Then $x^T Qx>0$ for every nonzero $x$ satisfying $x^\top Lx=0$ if and only if there exists $\varepsilon\in\R$ such that $Q+\varepsilon L\succ 0$.
\end{lemma}
In order to adapt the Lasserre--Netzer construction, we first need to make sure that the Hessian is positive definite on the whole space up to a perturbation. We write $\mathbb{D}$ for the symmetric tensor representation of the polarization of the determinant.
\begin{lemma}\label{lem: polynomial selection}
Let $\delta,R>0$, and let $p$ be a rank-one convex polynomial on $Q_R\subset \R^{2\times 2}$. Then there is a scalar polynomial $\Lambda\in \mathcal{P}(\R^{4})$ such that for all $A\in Q_R$ there holds
\begin{equation}\label{eq: polynomial perturbation}
 M(A)=\nabla^2p(A)+\delta \mathbb I+\Lambda(A)\mathbb{D} \succ 0. 
\end{equation}
\end{lemma}
\begin{proof}
    A matrix $H\in \R^{2\times 2}$ has rank at most one if and only if $\det H=0$.
In other words, for every $A$, by Lemma \ref{lem: Finsler's lemma}, there is a real number $\Lambda_A$ such that 
$$D^2p(A)+\frac{\delta}{2} \mathbb I+  \Lambda_A \mathbb{D}\succ 0.$$
For every $A\in \R^{2\times 2}$, there is an $r_A>0$ such that for all $B\in B(A,r_A)$, there holds
$$D^2p(B)+\delta \mathbb I+  \Lambda_A \mathbb{D}\succ \frac{\delta}{4}\mathbb{I}. $$
The balls $B(A,r_A)_{A\in Q_R}$ form an open cover of $Q_R$, whence we can choose a finite subcover $B_i:=B(A_i,r_{A_i})$, $i=1,...,l$, $l\in \N$. For this finite cover, we choose a partition of unity
$$\eta_i\in C^0(\R^{2\times 2},[0,1]), \quad \eta_i\equiv 0 \text{ on } B_i^c, \quad \sum_{i=1}^l\eta_i\equiv 1 \text{ on }Q_R,$$
and define a continuous function $\tilde{\Lambda}:=\sum_{i=1}^l\eta_i\Lambda_{A_i}.$
    By Weierstrass's theorem, we may choose a polynomial $\Lambda$ such that 
    $$\Vert \Lambda -\tilde{\Lambda}\Vert_{C^0(Q_R)}< \frac{\delta}{4},$$
    and, by $\vert \mathbb{D}\vert_{\mathrm{op}}=\frac{1}{2}$,
    $$D^2p+\delta \mathbb I+  \Lambda \mathbb{D}\succ \frac{\delta}{8}\mathbb I$$
    holds uniformly in $Q_R$. 
\end{proof}
Now that we have a positive definite matrix polynomial throughout $Q_R$, we can apply Prop.~\ref{Prop: Density of SOS polynomials}.
\begin{lemma}\label{lem: matrix Lasserre} Let $N\in \N$.
Let $M\in\mathcal{P}(\R^{m})^{N\times N}$ be a symmetric matrix polynomial that is pointwise positive definite on $Q_R$. For every $\varepsilon>0$, there exists $r_0$ such that, for every integer $r\ge r_0$,
$$M(u)+\varepsilon\Theta_{r,R}(u)\mathbb I_N$$
is an $\mathrm{SOS}$ matrix.
\end{lemma}
\begin{proof}
    We may first choose $\mu>0$ such that $M(u)-\mu\mathbb I_N\succ 0$ for all $u\in Q_R$. These matrices admit a positive definite square root matrix $C(u)$, which is continuous in $u$. By the Stone-Weierstrass theorem, there is a matrix polynomial $B\in \mathcal{P}(\R^{m})^{N\times N}$ approximating $C$, such that 
    $$\sup_{u\in Q_R}\Vert M(u)-B^T(u)B(u)-\mu\mathbb I_N\Vert_{\max} < \frac{\mu}{4N}.$$
    Let us write  $E:=M-B^TB$. Then this matrix is diagonally dominant, and we decompose $E$ using
    $$f_i=E_{ii}-(N-1)\frac{\mu}{2N}, \quad f_{ij}^+=\frac{\mu}{4N}+\frac{E_{ij}}{2}, \quad  f_{ij}^-=\frac{\mu}{4N}-\frac{E_{ij}}{2} $$
    into 
    \begin{align*}
        E= \sum_{i=1}^N f_i e_i\otimes e_i+\sum_{i<j}\biggl[
 f_{ij}^{+}(e_i+e_j)\otimes (e_i+e_j)
 +f_{ij}^{-}(e_i-e_j)\otimes (e_i-e_j)\biggr].
    \end{align*}
    This weighted split along rank-one matrices has the advantage that the scalar polynomials are nonnegative:
    \begin{align*}
         f_{ij}^\pm &\ge \frac{1}{2} \biggl(\frac{\mu}{2N}-\frac{\mu}{4N}\biggr)=\frac{\mu}{8N}>0,\\
         f_i &\ge \mu -\frac{\mu}{4N} -(N-1)\frac{\mu}{2N}=\mu \frac{(4N-1-2N+2)}{4N}>0.
    \end{align*}
    Set $\tau=\tfrac{\varepsilon}{2N-1}$ and apply Prop.~\ref{Prop: Density of SOS polynomials} to every $f$ in the decomposition, i.e., for every $f$ there is an integer $r_f\in \N$ such that for all $r\ge r_f$
    $$f+\tau \Theta_{r,R} \quad \text{is SOS}.$$
    Set $r_0:=\max_{f}\{r_f\}$, then for all $r\ge r_0$, writing $f_{i,\tau}=f_i+\tau\Theta_{r,R}$, $
f_{ij,\tau}^{\pm}=f_{ij}^{\pm}+\tau\Theta_{r,R}$,
    \begin{align*}
        &M+\varepsilon\Theta_{r,R}\mathbb I_N=M+\tau\Theta_{r,R}[1+2(N-1)]\mathbb I_N\\
        &=M + \tau\Theta_{r,R}\biggl(\sum_{i=1}^N e_i\otimes e_i+
 \sum_{i<j}[(e_i+e_j)\otimes (e_i+e_j)+(e_i-e_j)\otimes (e_i-e_j)]\biggr)\\
 &=B^TB+ \sum_{i=1}^N f_{i,\tau} e_i\otimes e_i \\
 & \qquad +\sum_{i<j}\biggl[
 f_{ij,\tau}^{+}(e_i+e_j)\otimes (e_i+e_j)
 +f_{ij,\tau}^{-}(e_i-e_j)\otimes (e_i-e_j)\biggr],
    \end{align*}which is a sum-of-squares polynomial. 
\end{proof}

We must also ensure that the perturbation is a Hessian matrix, which is not the case for the original construction. However, we can identify a polynomial that comes from a Hessian, such that the difference from $\varepsilon\Theta_{r,R}\mathbb I_N$ is a sum-of-squares polynomial. 

\begin{lemma}\label{lem: Hessian SOS}
For every integer $r\ge 1$ and every $u\in \R^N$, set
\begin{equation}\label{eq:psi}
 \psi_{r,R}(u)=\frac{1}{2}\vert u\vert ^2+
 \frac{R^2}{2(r+1)}\biggl(\frac{\vert u\vert^2}{R^2}\biggr)^{r+1}.
\end{equation}
 Then the polynomial matrix $ D^2\psi_{r,R}-\Theta_{r,R}\mathbb I_N$ is an $\mathrm{SOS}$ matrix polynomial.
\end{lemma}
\begin{proof}
We compute the second derivative of $\psi_{r,R}$
$$\nabla\psi_{r,R}=u+\bigl(\tfrac{\vert u\vert^2}{R^2}\bigr)^{r}u,
 \quad
 \nabla^2\psi_{r,R}=\mathbb I_N+\bigl(\tfrac{\vert u\vert^2}{R^2}\bigr)^{r}\mathbb I_N+
 \frac{2r}{R^2}\bigl(\tfrac{\vert u\vert^2}{R^2}\bigr)^{r-1} u\otimes u,$$
and infer
\begin{align}
    D^2 \psi_{r,R}-\Theta_{r,R}\mathbb I_N&= \biggl(\frac{\vert u\vert^{2r}}{R^{2r}}-\sum_{i=1}^N\frac{u_i^{2r}}{R^{2r}}\biggr)\mathbb I_N + \frac{2r}{R^2}\bigl(\tfrac{\vert u\vert^2}{R^2}\bigr)^{r-1} u\otimes u. \label{eq: Difference to Hessian} 
\end{align}
The multinomial theorem gives
$$ \frac{\vert u\vert^{2r}}{R^{2r}}-\sum_{i=1}^N\frac{u_i^{2r}}{R^{2r}}
 =\sum_{\substack{\vert \alpha\vert =r\\ \alpha\ne re_i\ (i=1,...,N)}}
   \binom r\alpha\biggl(\frac{u^\alpha}{R^r}\biggr)^2,$$
which is a sum-of-squares polynomial matrix if $r\ge 2$ (and zero otherwise). Thus, the first summand in \eqref{eq: Difference to Hessian} is an SOS matrix polynomial. Applying the multinomial theorem to the second summand gives
\begin{align*}
    \frac{2r}{R^2}\bigl(\tfrac{\vert u\vert^2}{R^2}\bigr)^{r-1} u\otimes u=\sum_{\vert \alpha\vert =r-1}\frac{2r}{R^{2r}}\binom{r-1}{\alpha}
   (u^\alpha u)\otimes (u^\alpha u),
\end{align*}which is again an SOS matrix polynomial. 
\end{proof}

We can now complete the proof of our approximation theorem. 
\begin{proof}[Proof of Thm.~\ref{Thm: Density of SOSRC}]
Let $f\colon \R^{m\times n}\to \R$ be continuous and rank-one convex. Let $\min\{m,n\}=2$ and $\varepsilon>0$. Let $K$ be a compact set and choose $R>0$ such that $K\Subset Q_{R/2}$. We first mollify $f$ with a standard compactly supported mollifier such that the result $f_\sigma$, $\sigma>0$, satisfies
$$\Vert f_\sigma-f\Vert_{C^0( Q_{2R})}<\frac{\varepsilon}{4}.$$
Mollification preserves rank-one convexity, since it is an average of translates of a rank-one convex function. Then, by the Weierstrass theorem, choose a polynomial $p$ such that
$$\biggl\Vert f_\sigma+\frac{\varepsilon}{8}\frac{1}{R^2}\vert \cdot \vert^2-p\biggl\Vert_{C^2( Q_{R})}<\min\biggl\{ \frac{\varepsilon}{8}, \frac{\varepsilon}{8mnR^2} \biggr\}.$$

In particular, the polynomial $p$ is rank-one convex on $Q_{R}$. By lemma \ref{Lem: compact open topology}, it suffices to find an SOSRC polynomial $\hat{p}$ such that
$$\sup_{A\in Q_{R/2}}\vert p(A)-\hat{p}(A)\vert<\frac{\varepsilon}{2}.$$
First, assume that $m=n=2$. Choose $\delta,\eta >0$ such that
$$\tfrac{1}{2}R^2\delta +\eta R^2< \frac{\varepsilon}{2}.$$
   For sufficiently large $r$, there exist a scalar polynomial $\Lambda$ and a polynomial matrix $B$ such that 
   $$D^2p+\delta \mathbb I+\Lambda \mathbb{D}+\eta  \Theta_{r,R}\mathbb I=B^TB$$
   by lemmas \ref{lem: polynomial selection} and \ref{lem: matrix Lasserre}. Moreover, Lemma \ref{lem: Hessian SOS} provides a polynomial matrix such that 
   $$D^2\psi_{r,R}-\Theta_{r,R}\mathbb I=C^TC.$$
Set
$$\hat{p}(A)=p(A)+\frac{\delta}{2}\Vert A\Vert^2+\eta \psi_{r,R}(A).$$
Then we have
\begin{align*}
    D^2\hat{p}+\Lambda \mathbb{D}&= D^2p+\delta\mathbb I+\eta D^2 \psi_{r,R}+\Lambda \mathbb{D} \\
    &=  D^2p+\delta\mathbb I+\Lambda \mathbb{D} +\eta \Theta_{r,R}\mathbb I+\eta D^2 \psi_{r,R}-  \eta  \Theta_{r,R}\mathbb I\\
    &= B^TB+\eta C^TC,
\end{align*}
and by our choice of the parameters, we obtain
\begin{align*}
    \sup_{A\in Q_{R/2}}\vert \hat{p}(A)-p(A)\vert&\le \frac{\delta}{2}4(R/2)^2+\eta  \frac{R^2}{2}+\eta\frac{R^2}{2(r+1)}\\
    &\le \frac{\delta}{2}4(R/2)^2+\eta R^2 <\frac{\varepsilon}{2}. 
\end{align*}
In summary, we have, for every $A\in Q_{R/2}$,
\begin{align*}
    \vert f(A)-\hat{p}(A)\vert &\le \vert f(A)-f_\sigma(A)\vert +\vert f_\sigma(A) - p(A)\vert +\vert p(A)-\hat{p}(A)\vert\\
    &\le \frac{\varepsilon}{4} +\frac{\varepsilon}{8}+\frac{\varepsilon}{8}\frac{1}{R^2}\vert A\vert^2+ \frac{\varepsilon}{2}\le \varepsilon.
\end{align*}
    Finally, assume only that $\min\{m,n\}=2$. Then the proof proceeds analogously, except for lemma \ref{lem: polynomial selection}. The reason is that Finsler's lemma is generally not compatible with multiple constraints $G$. However, with either two columns or two rows, we can alternatively argue with Calderón's Theorem \cite{calderon1973note}, which states that, for every $A$, the Hessian $D^2p(A)[a\otimes b,a\otimes b]$ is a sum of squares in $(a,b)$. This allows us to again add a positive-definite quadratic form given by $\varepsilon \mathbb I$, and extend to a positive definite form on the full-rank matrices. We then obtain a continuous representative by a partition of unity and approximate its coefficient functions by polynomials using the Weierstrass theorem, as before. The difference between this extension and $D^2p(A)[H,H]$ is then of the form $\Lambda:M_2(H)$, where $M_2$ is the vector of second-order minors and $\Lambda\in \mathcal{P}(\R^{mn})^{\binom{\max\{m,n\}}{2}}$ is a multiplier. The proof continues analogously to the $m=n=2$ case, by choosing $Q_R$ such that $K\Subset Q_{R/\sqrt{2\max\{m,n\}}}$.
\end{proof}

Finally, let us discuss how density fails if $m,n\ge 3$. In a nutshell, Calderón's theorem is not applicable anymore in these directions. More precisely, there is  a positive semidefinite biquadratic form $F$ that cannot be expressed as the sum of squares of bilinear forms \cite{choi1975positive}. Choose a quadratic polynomial $q$ on $\R^{3\times3}$ such that $D^2q[a\otimes b,a\otimes b]=F(a,b)$. This $q$ is rank-one convex. Now, suppose there is a sequence of SOSRC polynomials $(q_j)_{j\in \N}$ approximating $q$ locally uniformly. Take a nonnegative smooth mollifier $\varphi$ with compact support and define
$$F_j(a,b):=\int_{\R^{3\times 3}} \varphi(X) 
D^2q_j(X)[a\otimes b,a\otimes b]\dd X.$$
Each $F_j$ is s sum of squares of bilinear forms. Integration by parts twice shows that $F_j$ converges coefficientwise to $F$. Since the fixed-degree SOS cone is closed, $F$ must be SOS, which is a contradiction. This example extends naturally to higher dimensions

\section{Representation of SOSRC polynomials}
\label{Sec: Generation}

We begin with a simple lemma on the form of the Hessian of SOSRC polynomials.
\begin{lemma}\label{lem: Hessian SOS RC}
    Let $p\in \mathcal{P}(\R^4)$ be an SOSRC polynomial of degree $2d$. Then there exist an integer $N\in \N$, a polynomial $\Lambda\in \mathcal{P}(\R^4)$ of degree at most $2d-2$, and polynomials
    $v_{\nu,i,j}\in \mathcal{P}(\R^4)$ of degree $\deg v_{\nu,i,j}\le d-1 $ for $\nu=1,...,N$, $i,j\in \{1,2\}$ such that
    \begin{equation}\label{eq: Hessian Decomposition}
        D^2p(A)[H,H]=\sum_{\nu=1}^N\biggl(\sum_{i,j=1}^2v_{\nu,i,j}(A)H_{ij}\biggr)^2+\Lambda(A)\det H
    \end{equation}
    for all $A,H\in \R^{2\times 2}\cong \R^4$. 
\end{lemma}
\begin{proof}
    By the SOS property, there exist an integer $N\in \N$ and polynomials $q_\nu(A,a,b)$ such that
    $$D^2p(A)[a\otimes b,a\otimes b]=\sum_{\nu=1}^Nq_\nu (A,a,b)^2.$$
Since the highest-degree terms in a sum of squares cannot cancel, the separate homogeneity of the left-hand side implies that every $q_\nu$ is bihomogeneous of degree one in $a$ and degree one in $b$.
    
  We collect the coefficients of the linear forms into polynomial matrices $v_{\nu}(A)$. We can simply extend the resulting right-hand side, replacing $a_ib_j$ with arbitrary matrix components $H_{ij}$, and thereby obtain the first summand in \eqref{eq: Hessian Decomposition}. 
    The difference between the left-hand side and this first summand vanishes on rank-one matrices. By \cite[Cor.~2]{bandeira2016role} there is a scalar  multiplier, which can then be seen to be a polynomial of
degree at most $2d-2$ by coefficient matching.
\end{proof}
An immediate observation is that $p$ is convex if $\Lambda\equiv 0$, and polyconvex if $\Lambda$ is constant. The converse implications do not hold in general. Nevertheless, these sufficient conditions will allow us to identify and discard certain polyconvex functions more easily in the numerical scheme.

\begin{proof}[Proof of Thm.~\ref{Thm: Representation}]
So far, we only characterized the Hessian of $p$. We would now like to obtain a direct polynomial characterization of $p$ itself. Any such characterization must incorporate the corresponding compatibility conditions. Let us identify $\R^{2\times 2}\cong \R^4$, and write
$$\mathcal{Q}(A):=\sum_{\nu=1}^N v_\nu(A)\otimes v_\nu(A)+\Lambda(A)\mathbb{D}\in \R^{4\times 4}.$$
Then for all $1\le i,j,k\le 4$, we must have
$$\partial_k\mathcal{Q}_{ij}= \partial_i \mathcal{Q}_{kj}$$
in order for the right-hand side to be a Hessian. 
Let us expand
$$v_\nu(A)=\sum_{\vert\alpha\vert\le d-1}v_{\nu,\alpha}A^\alpha, \ v_{\nu,\alpha}\in \R^4,
 \qquad \Lambda(A)=\sum_{\vert\beta\vert\le 2d-2}\lambda_\beta A^\beta, \ \lambda_\beta \in \R$$
and group the terms according to the respective monomials
$$ Q_\gamma:=\sum_{\nu=1}^N\sum_{\alpha+\beta=\gamma}
 v_{\nu,\alpha}\otimes v_{\nu,\beta}+\lambda_\gamma \mathbb{D}.$$
Then the compatibility condition for the coefficients takes the form 
$$(\gamma_k+1)(Q_{\gamma+e_k})_{ij}  =(\gamma_i+1 )(Q_{\gamma+e_i})_{kj}$$
for all $\vert \gamma\vert \le 2d-3$ and $i,j,k\in \{1,...,4\}$.
A family of polynomials $v_\nu,\Lambda$ satisfying these compatibility conditions can be integrated to obtain
$$p(X)=a+b\cdot X+
 \sum_{\vert \beta\vert \le2d-2}
 \frac{X^TQ_\beta X}{(\vert \beta\vert +1)(\vert \beta\vert +2)}X^\beta$$
 where $a\in\R$ and $b\in\R^4$ are arbitrary.
Regrouping these terms by degree yields \eqref{eq: integrated p}.
It remains to show that $p_{\mathrm{sos}}$ is actually a sum of squares. Consider the $d\times d$ matrix
$$\mathbb{W}_{i+1,j+1}=\frac{1}{(i+j+1)(i+j+2)}=\int_0^1 (1-t)t^{i+j}\dd t. $$
Thus for every $A\in \R^d$ there holds
\begin{align*}
    A^T\mathbb{W}A &=\int_0^1 (1-t) \sum_{i,j=0}^{d-1}A_{i+1}A_{j+1}t^{i+j}\dd t\\
    &= \int_0^1 (1-t) \biggl( \sum_{i=0}^{d-1}A_{i+1}t^{i}\biggr)^2\dd t\ge 0,
\end{align*}
and therefore there is a square root $\mathbb{W}=\mathbb{V}^T\mathbb{V}$. Hence,
\begin{align*}
 p_{\mathrm{sos}}(X)=\sum_{\nu=1}^N\sum_{j=1}^d
 \biggl(\sum_{i=0}^{d-1}\mathbb{V}_{j,i+1}X^Tv^{\mathrm{hom},i}_{\nu}(X)\biggr)^2.
\end{align*}
\end{proof}

\begin{proof}[Proof of Cor.~\ref{Cor: Generation of SOSRC functions}]
For each $\nu$, arrange the coefficient vectors $v_{\nu,\alpha}\in\R^4$ in the order used in $m_{d-1}$, and denote the resulting vector by $w_\nu\in\R^{N_d}$. Then $w_\nu^Tz_d(X,H)=v_\nu(X)^TH$.
Setting $\G=\sum_\nu w_\nu w_\nu^T\succeq 0$, the Hessian decomposition gives
$$D^2p(X)[H,H]=z_d(X,H)^T\mathbb Gz_d(X,H)+\Lambda(X)\det H.$$
Omitting the affine part of $p$ yields $\mathscr{F}_d(C,\mathbb G,\lambda)\equiv0$. Conversely, for any feasible triple, take $\mathbb G=\sum_\nu w_\nu w_\nu^T$ and substitute $H=a\otimes b$. Since $\det(a\otimes b)=0$, the resulting identity yields the SOSRC property for $p_C+p_{\mathrm{aff}}$.
\end{proof}

\section{An SDP approach to Morrey's problem}
\label{Sec: Morrey SDP}

In this section, we first explain how to use Cor.~\ref{Cor: Generation of SOSRC functions} to reformulate Morrey's problem as a semidefinite optimization problem at any choice of base point $A\in \R^{2\times 2}$, test function $u\in C^\infty(\mathbb{T}^2,\R^2)$, and fixed degree $d\in \N$. We then derive the dual problem of \eqref{eq: semidefinite program}, prove strong duality, and investigate strict feasibility of the dual problem. 
Thereafter, we present a numerical scheme for the SDP, and finally discuss the results of our test run.

\subsection{Morrey's problem as an SDP}

First, coefficient matching, normalization of $\Lambda$ and the objective function $L_u$ are linear in the input vector $C$. The constraint $\G\succeq 0$ is a semidefinite constraint. Hence \eqref{eq: semidefinite program} is indeed a semidefinite program. 
We first prove a more direct version of the SDP.
\begin{lemma}\label{Lem: general version SDP}
     For $u\in C^\infty(\mathbb{T}^2,\R^2)$ and an $\mathrm{SOSRC}$ polynomial $p_C$ of degree $2d$ with coefficient vector $C$, we introduce the objective function $$L_u(C)=\sum_{\vert\alpha\vert=2}^{2d}C_\alpha\int_{\mathbb T^2}(Du)^\alpha\dd x,$$ and the semidefinite program
    \begin{align}\label{eq: general version semidefinite program}
        \begin{aligned}
           \varrho_d&(u)=\inf_{C,\G,\lambda} L_u(C)\\
            &\mathrm{subject\ to} \ \mathscr{F}_d(C,\G,\lambda)=0, \ \operatorname{tr}\G\le 1,\ \mathrm{and}\ \G\succeq 0.
        \end{aligned}
    \end{align}
    Then there is a rank-one convex function on $\R^{2\times 2}$ failing to be quasiconvex at $0$ with witness $u$ if and only if there is a $d\in \N$ such that $\varrho_d(u)<0$.
\end{lemma}
\begin{proof}
Let us first assume that there is a $u$ and a degree $d$ such that $\varepsilon:=-\varrho_d(u)>0$. Then there are compatible $\G^*,\Lambda^*,C^*$ such that
$$L_u(C^*)=\int_{\mathbb{T}^2} p_{C^*}(Du)\dd x< -\frac{\varepsilon}{2}. $$
This immediately implies that $\tr \G^*>0$, since otherwise $\G^*=0$, and then the compatibility conditions force $\Lambda$ to be constant. Thus $p_C$ is a multiple of the determinant, up to affine terms, hence a null Lagrangian, and therefore quasiconvex and $L_u(C^*)=0$.

For the converse direction, assume that there exists a rank-one convex function $f$ that fails quasiconvexity at $0$ with witness $u$, i.e., 
$$\int_{\mathbb{T}^2} f(Du)-f(0)\dd x= -\varepsilon<0.$$
W.l.o.g., we may assume that $f(0)=0$ (otherwise replace $f$ by $f-f(0)$). Then, by Thm.~\ref{Thm: Density of SOSRC}, there is an SOSRC polynomial $p\in \mathcal{P}(\R^{2\times 2})$ of degree $2d$ such that
$$\int_{\mathbb{T}^2} p(Du)-p(0) \dd x<-\frac{\varepsilon}{2}.$$
Again, we may take $p$ such that $p(0)=0$. In particular, for the non-affine part $p_C$ of $p$, there holds
$$L_u(C)=\int_{\mathbb{T}^2} p_C(Du) \dd x<-\frac{\varepsilon}{2}.$$
By Thm.~\ref{Thm: Representation} and Cor.~\ref{Cor: Generation of SOSRC functions}, there are matching $\G,\Lambda$ such that $\mathscr{F}_d(C,\G,\lambda)=0$, and $\G\succeq 0$. Finally, if $\tr \G>1$, we can rescale $p$ by $\frac{1}{\tr \G}$, yielding
$$\varrho_d(u)\le L_u(\tfrac{1}{\tr \G}C)=\int_{\mathbb{T}^2} p_{\frac{1}{\tr \G}C}(Du) \dd x<-\frac{\varepsilon}{2\tr \G}<0.$$
\end{proof}

Note that if $p$ violates the inequality at a base point $A\ne 0$, the translated
polynomial $q(X) = p(A+X)$ is still SOSRC of degree at most $2d$. Removing its affine part
therefore gives a negative feasible point in \eqref{eq: general version semidefinite program} for the same $u$.

Let us now continue with the main proof.
\begin{proof}[Proof of Thm.~\ref{Thm: Morrey Test}]
We show that the SDP in \eqref{eq: general version semidefinite program} is equivalent to \eqref{eq: semidefinite program}.
For 
$$ p_C(X)= \int_0^1(1-t) \biggl[ z_d(tX,X)^T\G z_d(tX,X) +\Lambda(tX)\det X \biggr]dt,$$
it holds that $$L_u(C)=\langle \mathbb{M}_u,\G\rangle+b_u^T\lambda.$$
Note that a feasible $C$ exists if and only if $\G$ and $\Lambda$ are compatible. Moreover, the constant part of $\Lambda$ multiplied by the determinant is a null Lagrangian, and thus does not alter $L_u$ at all, and we can discard it like the affine part of $p$. 
Since we now normalize $b_u^T\lambda=-1$, we have $\varrho_d(u)<0$ if and only if $\rho_d(u)<1$. 
\end{proof}
The advantage of  \eqref{eq: semidefinite program} over \eqref{eq: general version semidefinite program} is that it already excludes admissible points that cannot violate the quasiconvexity inequality. In particular, we discard not just the null Lagrangians, but also SOSRC polynomials that can never witness the failure of Jensen's inequality for $u$. More precisely, the quadratic part of $L$ is always nonnegative, since the matrix $\mathbb{M}_u$ is positive semidefinite. To see this, let $A\in \R^{N_d}$. Then
\begin{align*}
A^T\mathbb{M}_uA&=\int_0^1(1-t)\int_{{\mathbb{T}}^2}\sum_{i,j=1}^{N_d} z_{d,i}(tDu,Du)z_{d,j}(tDu,Du)A_iA_j\dd x\dd t\\
& =\int_0^1(1-t)\int_{{\mathbb{T}}^2} (z_d(tDu,Du)^TA)^2\dd x\dd t\ge 0.
\end{align*}
Hence, only the nonconstant part of $\Lambda$ can potentially cause failure of quasiconvexity. We therefore can limit our investigation to those coefficients $\lambda $ that admit a negative $b_u^T\lambda$. 

\subsection{Duality in the Morrey SDP}

We begin by introducing the following notation. We define linear operators

\begin{align}\label{eq: Linear Operators}
   \begin{aligned}
        \mathcal{H}&\colon \mathcal{P}_{2d}\to \Sym, \ p\mapsto D^2 p,\\
    \mathcal{C}&\colon \R^{{N_d}\times {N_d}}_{\mathrm{Sym}}\to \Sym, \ \G\mapsto (m_{d-1}(X)\otimes \mathbb I_4)^T\G(m_{d-1}(X)\otimes \mathbb  I_4),\\
    \mathcal{D}&\colon \{\Lambda\in\mathcal{P}_{2d-2}\mid \Lambda(0)=0\} \to \Sym, \ \Lambda \mapsto \Lambda \mathbb{D}.
   \end{aligned}
\end{align}
Then we can write the feasibility operator as
$$\mathscr{F}_d(C,\G,\lambda)= \mathcal{H}(p_C)- \mathcal{C}(\G)-\mathcal{D}(\lambda). $$
We can equip $\Sym$ naturally with an inner product by
$$\biggl\langle \sum_{\vert \alpha\vert \le 2d-2}S_\alpha X^\alpha, \sum_{\vert \beta\vert \le 2d-2}R_\beta X^\beta\biggr\rangle:=\sum_{\vert \alpha\vert \le 2d-2}\tr(S_\alpha^T R_\alpha).$$
Now, we derive the form of the dual problem.

\begin{lemma}
    The dual problem of \eqref{eq: semidefinite program}
    is given by
  \begin{align}\label{eq: dual program} \tag{DSP}
        \begin{aligned}
          &\sup_{Y,\eta}\, \eta\\
           &\mathrm{subject\ to} \  \mathcal{H}^*(Y)=0, \ \mathcal{D}^*(Y)=\eta b_u, \ \mathrm{and}\  \mathbb{M}_u-\mathcal{C}^*(Y)\succeq 0.
        \end{aligned}
    \end{align}
\end{lemma}
\begin{proof}
Let us formulate the primal problem with the feasibility operator: 
\begin{align*}
            &\rho_d(u)=\inf_{C,\G,\lambda}\, \langle \mathbb{M}_u,\G\rangle \\
            &\mathrm{subject\ to} \ \mathscr{F}_d(C,\G,\lambda)=0, \ b_u^T\lambda= -1, \ \mathrm{and}\ \G\succeq 0.
        \end{align*}
    We use the Lagrangian method and introduce the dual variables $\eta\in \R$ and $Y\in \Sym$. The Lagrangian is given by
    $$\mathscr{L}(C,\G,\Lambda;\eta,Y):=\langle \mathbb{M}_u,\G\rangle+\langle Y,\mathcal{H}(p_C)- \mathcal{C}(\G)-\mathcal{D}(\Lambda)\rangle +\eta(b_u^T\lambda+1).$$
Varying over $C$, we infer $\mathcal{H}^*(Y)=0$.
Next, varying over $\Lambda\in \mathcal{P}_{2d-2}$, $\Lambda(0)=0$, we deduce
$$\mathcal{D}^*(Y)=\eta b_u.$$
Finally, varying over $\G\succeq 0$, we get 
$$\langle \mathbb{M}_u-\mathcal{C}^*(Y),\G\rangle\ge 0.$$
This forces $\mathbb{M}_u-\mathcal{C}^*(Y)\succeq 0$. 
\end{proof}
We remark that the adjoint operators have an explicit form. For $Y(X)=\sum_{\vert \gamma\vert \le 2d-2}Y_\gamma X^\gamma$, they are given by (terms with a negative component in
$\nu-e_i-e_j$ are omitted)
\begin{align*}
    (\mathcal{C}^*(Y))_{\alpha\beta}
=Y_{\alpha+\beta}, \qquad 
&(\mathcal{D}^*Y)_\gamma
=\langle Y_\gamma,\mathbb{D}\rangle
=(Y_\gamma)_{14}-(Y_\gamma)_{23}, \\
(\mathcal{H}^*Y)_\nu
= &\sum_{i,j=1}^{4} \nu_i(\nu_j-\delta_{ij}) (Y_{\nu-e_i-e_j})_{ij}.
\end{align*}

Next, we show that the primal problem contains an interior point. Consequently, the supremum in the dual problem is attained, and we have strong duality, meaning that the values of the primal and dual problem coincide. 
\begin{lemma}
    The primal problem \eqref{eq: semidefinite program} contains an interior point if $b_u\ne 0$.
\end{lemma}
\begin{proof}
    \textit{1.~Step.~}$\mathcal{C}$ \textit{is surjective: }For the following arguments, it is convenient to think of $\mathbb{F}\in \R^{{N_d}\times {N_d}}_{\mathrm{Sym}}$ as a matrix in $\mathbb{S}^{{N_d/4}\times {N_d/4}}_{\mathrm{Sym}}$, the entries of which are in $\R^{4\times 4}$. Similarly, one should think of $m_{d-1}\otimes \mathbb  I_4$ as a $N_d/4$-dimensional vector with entries $x^\alpha  \mathbb I_4$. Let $R\in \Sym$. For any coefficient matrix $R_\gamma$, choose multi-indices with $\vert \alpha\vert,\vert\beta\vert\le d-1$ and $\alpha+\beta=\gamma$. If $\alpha=\beta$, set
    $\mathbb{F}_{\alpha\alpha}=R_\gamma$, otherwise $\mathbb{F}_{\alpha\beta}=\mathbb{F}_{\beta\alpha}=\tfrac{1}{2}R_\gamma$. Then the coefficient of $X^\gamma$ in $\mathcal{C}(\mathbb{F})$ is $R_\gamma$.\\
\textit{2.~Step.~A simple perturbation: }Let $p_0=(1+\vert X\vert^2)^d-1$. Its Hessian is given by
$$\mathcal{H}(p_0)=2d(1+\vert X\vert^2)^{d-1} \mathbb I_4+4d(d-1)(1+\vert X\vert^2)^{d-2}X\otimes X. $$
The matrix $X\otimes X$ is positive semidefinite, therefore the second summand is positive semidefinite. The first summand is positive definite. This can be seen by employing the multinomial theorem
$$(1+\vert X\vert^2)^{d-1}=\sum_{\vert\alpha\vert\le d-1}c_\alpha X^{2\alpha}, \qquad c_\alpha=\binom{d-1}{\vert \alpha\vert}\frac{\vert \alpha\vert!}{\vert \alpha_1\vert!\vert \alpha_2\vert!\vert \alpha_3\vert!\vert \alpha_4\vert!}>0.$$
Writing $C=(c_\alpha)_\alpha$, we see that 
$$2d(1+\vert X\vert^2)^{d-1} \mathbb I_4=2d(m_{d-1}(X)\otimes \mathbb I_4)^T(\diag(C)\otimes  \mathbb I_4)(m_{d-1}(X)\otimes \mathbb I_4),$$
and the matrix $\diag(C)\otimes \mathbb  I_4$ is positive definite, which is immediate if we again write it as a block-diagonal matrix with $N_d/4$ blocks of the form $c_\alpha  \mathbb I_4$. In particular, there is a positive definite $\G_d$ such that $\mathcal{C}(\G_d)=\mathcal{H}(p_0)$. \\
\textit{3.~Step.~Applying Weyl's inequality: }Since $b_u\ne 0$, while $(b_u)_0=0$, we can find an admissible $\lambda$ by setting $\lambda=-\frac{b_u}{\Vert b_u\Vert^2}$. By step 1, we find a matrix $\mathbb{F}$ such that $\mathcal{C}(\mathbb{F})=-\Lambda \mathbb{D}$. Setting $\G_t:=t\G_d+\mathbb{F}$, we get
$$\mathcal{C}(\G_t)=\mathcal{H}(tp_0)-\mathcal{D}(\Lambda),$$
thus $\G_t,t p_0,\Lambda$ are feasible. Finally, we can use Weyl's spectral inequality (see, e.g.~\cite[Cor.~III.2.6]{Bhatia1997MatrixAnalysis}) to estimate the minimal eigenvalue $\mu_{\min}$ as
$$\mu_{\min} (\G_t)\ge t\mu_{\min}(\G_d)-\vert \mathbb{F}\vert_{\mathrm{op}}.$$
Hence, $\G_t$ is positive definite for sufficiently large $t$ by step 2.
\end{proof}
The existence of a strictly feasible point for the dual problem turns out to be a rather complicated condition on the test function $u$, namely \eqref{eq: Strict duality condition}.
\begin{lemma}
    The dual problem \eqref{eq: dual program} contains an interior feasible point if and only if condition \eqref{eq: Strict duality condition} holds. 
\end{lemma}
\begin{proof}
We assume that $b_u\ne 0$, otherwise the SDP is trivial. 
Let us set
$$\mathcal{T}_u=\{Y\in \Sym\mid \mathcal{H}^*(Y)=0,\mathcal{D}^*(Y)\in \operatorname{span}(b_u)\}.$$
    The feasible set can then be represented by 
    $$\{(Y,\eta)\mid Y\in \mathcal{T}_u, \mathcal{D}^*(Y)=\eta b_u, \mathbb{M}_u+\mathcal{C}^*(-Y)\succeq 0\}.  $$
This set contains an interior point precisely if there is a $Y\in \mathcal{T}_u$ such that $\mathbb{M}_u+\mathcal{C}^*(-Y)\succ 0$. By \cite[Thm.~1]{balakrishnan2003semidefinite}, this is equivalent to the non-existence of a nonzero $\G\succeq 0$ with the properties $\mathcal{C}(\G)=\mathcal{H}(p_C)+\mathcal{D}(\Lambda)$, $b_u^T\lambda=0$ and $\langle \mathbb{M}_u,\G\rangle\le 0$. Since $\mathbb{M}_u$ is positive semidefinite, the last inequality can only be an actual equality.
\end{proof}
The natural question arising is whether there are any nontrivial test functions $u$ that do not admit the strict dual feasibility condition \eqref{eq: Strict duality condition}. 
In fact, there are many such examples, one of the simplest ones being 
$$u(x_1,x_2)=(\sin(2\pi x_1),\sin(2\pi x_2))^T.$$
Choose for instance the monomial $X_1X_4^{2d-3}$ for $d\ge 2$. Then
\begin{align*}
    &\int_{\mathbb{T}^2} \det(Du)\partial_1 u_1 (\partial_2 u_2)^{2d-3}\dd x\\
    &=(2\pi)^{2d}\biggl(\int_0^1 \cos(2\pi x_1)^2\dd x_1\biggr)\biggl(\int_0^1 \cos(2\pi x_2)^{2d-2}\dd x_2\biggr)>0.
\end{align*}
Let $e$ denote the coordinate indexed by $(\alpha,i)=((0,d-1,0,0),2)$, i.e., $e^Tz(x,h)=x_2^{d-1}h_2$. Set $G_0=2d(2d-1)ee^T$, and $p_0(x)=x_2^{2d}$. take $\lambda=0$, then
$$\mathcal{C}(G_0)=\mathcal{H}(p_0),
\qquad
b_u^T\lambda_0=0,$$
and $G_0\succeq 0$ is nonzero. However, $(Du)_{12}=0$ holds everywhere, so
\begin{align*}
\langle M_u,G_0\rangle
&= \int_{\mathbb T^2}(Du)_{12}^{2d}\dd x =0.
\end{align*}
Such a simple counterexample immediately raises the question of whether there are any $u$ such that an interior point in \eqref{eq: dual program} exists. Such functions are more difficult to identify, since they depend heavily on the chosen degree bound $2d$. 
One such example we have found for degrees $d=1,2,3$ is
$$u(x_1,x_2)=\bigl(\sin(2\pi x_1)+\sin(2\pi x_2),\,
\sin(2\pi(x_1+x_2))+\sin(2\pi(x_1-x_2))\bigr)^T.$$
The corresponding dual feasible points were obtained numerically and subsequently verified using rational arithmetic. We omit the rather lengthy certificates here, as these examples are only intended to illustrate that dual strict feasibility can indeed occur, but should be obtained through an additional numerical procedure before solving the actual SDP.

\subsection{A numerical scheme for Morrey's conjecture}

Let us first express \eqref{eq: semidefinite program} in a more general form. The compatibility equation reads $\mathcal{C}(\G)=\mathcal{H}(p_C)-\mathcal{D}(\Lambda)$. Let us join these conditions in a linear operator via $\mathcal{B}(C,\lambda):=\mathcal{H}(p_C)-\mathcal{D}(\Lambda)$, then we require $\mathcal{C}(\G)\in \operatorname{R}(\mathcal{B})$. Since $\operatorname{R}(\mathcal{B})=\ker(\mathcal{B}^*)^\perp$, we obtain $1,...,k:=\dim \ker(\mathcal{B}^*)$ linear equations
$$\langle \mathcal{C}^*Y_j,\G\rangle=0,$$
where $Y_1,...,Y_k$ form a basis of $\ker(\mathcal{B}^*)$. 
For the normalization, we seek $Y_0\in\Sym$ such that $Y_0\perp\operatorname{R}(\mathcal{H})$ and $Y_0-\mathcal{D}(b_u)\perp\operatorname{R}(\mathcal{D})$. Such a $Y_0$ exists because $\operatorname{R}(\mathcal H)\cap\operatorname{R}(\mathcal D)=\{0\}$, as the Hessian compatibility conditions force any multiplier $\Lambda$ satisfying $D^2p=\Lambda\mathbb{D}$ to be constant, whereas $\Lambda(0)=0$ then forces it to be zero. Moreover, $\mathcal D^*\mathcal D=\id$, since $\|\mathbb D\|_2^2=1$. Then
\begin{align*}
    \langle Y_0,\mathcal{C}(\G)\rangle &= \langle Y_0,\mathcal{H}(p_C)-\mathcal{D}(\Lambda)\rangle= \langle \mathcal{H}^*Y_0,p_C\rangle-\langle\mathcal{D}^*Y_0,\lambda\rangle=-b_u^T\lambda=1.
\end{align*}
Thus \eqref{eq: semidefinite program} can equivalently be written as
\begin{align}\label{eq: reduced SDP}
        \begin{aligned}
           &\rho_d(u)=\inf_{\G}\, \langle \mathbb{M}_u,\G\rangle \\
            &\mathrm{subject\ to} \  \langle \mathcal{C}^*Y_j,\G\rangle=0, j=1,...,k, \  \langle \mathcal{C}^*Y_0,\G\rangle= 1, \ \mathrm{and}\ \G\succeq 0.
        \end{aligned}
    \end{align}
These considerations lead us to the following numerical scheme.
 \begin{Scheme}\textcolor{white}{.}
    \begin{itemize}
    \item[In:] Test function $u$ (as Fourier series) and maximal degree $2d$
\item[1.] Compute the moment data $\mathbb M_u$ and $b_u$ from $Du$.
\item[2.] Eliminate the polynomial variables to obtain an equivalent
       semidefinite program in the matrix $\G$.
\item[3.] Solve the reduced SDP using a primal-dual interior-point method.
\item[4.] Reconstruct the polynomial $p_C$ from the latest computed $C$ and verify the resulting
       rank-one convexity and Jensen violation. 
        \item[Out:]  (Almost) optimal value $\rho_d(u)$ and a corresponding polynomial $p_C$ 
    \end{itemize}
\end{Scheme}

Before solving the semidefinite program, let us revisit the Slater conditions. While primal strict feasibility is essential and assumed to be satisfied (otherwise we can discard $u$ in the first place), dual strict feasibility is not a necessity in theory. The latter can still improve the numerical stability of the interior-point method, which is why we restrict the numerical search to test functions for which dual strict feasibility can be certified with a sufficiently positive numerical margin.
Such a condition can be imposed directly by solving another SDP to compute the dual Slater margin.

The reduced SDP is solved via a primal-dual interior-point method. Starting from positive-definite iterates for both primal and dual problems, the algorithm improves the feasibility and decreases the duality gap through Newton iterations. We use the predictor-corrector approach by Mehrotra \cite{mehrotra1992implementation}, where the prediction step is affine and the corrector step moves toward the central path. The resulting linear systems at each iteration are solved by computing the Schur complement of the primal-dual system. The procedure terminates if the duality gap is below a certain threshold. 

Concerning the computational cost of our algorithm, let $N_d$ denote the order of $\G$ and $N_e$ be the number of equality constraints. It holds that
$$  N_d=4\binom{d+3}{4}=\mathcal{O}(d^4),
    \qquad N_e=\mathcal{O}(d^4).$$
The Schur complement multiplication requires $\mathcal{O}(N_d^3+N_e N_s)=\mathcal{O}(d^{12})$ operations, where $N_s=10\binom{2d+2}{4}=\mathcal{O}(d^4)$ is the dimension of the coefficient space of $\Sym$. Let $N_i$ be the number of interior-point iterations and $N_{cg}$ the conjugate gradient iterations per Schur complement, then the resulting cost is
$$ \mathcal{O}(d^{12} +N_i(1+N_{cg})d^{12}).$$
In practice, the numbers $N_i,N_{cg}$ are of order $\mathcal{O}(1)$.
Consequently, we were only able to test-run our SDP solver for low degrees (up to $2d=12$) and could not obtain a counterexample in this way.

\section{Extensions to isotropic and symmetric rank-one convexity}
\label{Sec: Corollaries}

We can apply the reasoning behind our density and semidefinite programming approach to other variants of Morrey's problem as well. In particular, we first consider rank-one convex functions with additional invariances (such as isotropy), and second, we regard the symmetric Morrey problem, which is more related to geometrically linear elasticity. 

\subsection{Corollaries on isotropic integrands}

We can extend the density result simply to integrands while imposing an additional invariance by averaging with respect to the Haar measure. 

\begin{Proposition}\label{Prop: Invariances}
    Let $G$ be a compact subgroup of $\mathrm{GL}(\R^{2\times 2})$. Assume that the action preserves the rank-one cone and the class of $\mathrm{SOSRC}$
polynomials for every $g\in G$. Then every continuous, $G$-invariant rank-one convex function
$f\colon \R^{2\times 2}\to\R$ can be locally uniformly approximated by $G$-invariant $\mathrm{SOSRC}$ polynomials.
\end{Proposition}
\begin{proof}
Let $K\subset \R^{2\times 2}$ be compact and $\varepsilon>0$. Since $G$ is compact, the orbit $G\cdot K:=\{gA\mid g\in G, A\in K\}$ is compact. By Thm.~\ref{Thm: Density of SOSRC}, there exists an SOSRC polynomial $p_0$ such that
$$\sup_{A\in G\cdot K}\vert p_0(A)-f(A)\vert <\varepsilon.$$
Define
$$p(A):=\int_G p_0(gA)\,\mathrm d\mu_G(g),$$
where $\mu_G$ denotes the normalized Haar measure on $G$. Then $p$ is a polynomial and is $G$-invariant by invariance of the Haar measure. Moreover, $p$ is SOSRC, since averaging the
corresponding positive semidefinite matrices in the SOS part of $p$ (cf.~Thm.~\ref{Thm: Representation}) preserves positive semidefiniteness by assumption. Finally, by the $G$-invariance of $f$, for every $A\in \R^{2\times 2}$, it holds that
\begin{align*}
\vert p(A)-f(A)\vert &\le \int_G
\vert p_0(gA)-f(gA)\vert \dd \mu_G(g)<\varepsilon
\end{align*}
for every $A\in K$.
\end{proof}
\begin{corollary}\label{Cor: Rotaional invariances}
    Let $f\colon \R^{2\times 2}\to \R$ be rank-one convex. 
    \begin{itemize}
        \item[(i)] If $f$ is invariant under left (resp.~right) rotations, i.e., for all $R\in \mathrm{SO}(2)$, $A\in \R^{2\times 2}$ there holds $f(RA)=f(A)$ (resp.~$f(AR)=f(A)$), then $f$ can be approximated in the local uniform topology by $\mathrm{SOSRC}$ polynomials that are invariant under left (resp.~right) rotations. 
        \item[(ii)] If $f$ is isotropic, i.e., for all $R,S\in \mathrm{SO}(2)$, $A\in \R^{2\times 2}$ it holds $f(RAS^T)=f(A)$, then $f$ can be approximated in the local uniform topology by $\mathrm{SOSRC}$ polynomials that are isotropic. 
    \end{itemize}
\end{corollary}
\begin{proof}
    For (\emph{i}), let $G\subset \mathrm{GL}(\R^{2\times 2})$ be the group defined left (resp.~right) multiplications with rotations. This group is compact. The action preserves the rank of a matrix and the determinant exactly, since $R(a\otimes b)S^T=(Ra)\otimes(Sb)$. In particular, it preserves the SOSRC property of polynomials. Then apply Prop.~\ref{Prop: Invariances}.\\
    For (\emph{ii}), apply the same argument to the product group with action defined $(R,S)\mapsto RAS^T$.
\end{proof}
We now focus on the isotropic case. Here, we establish a simplified representation theorem. 
\begin{lemma}
    Let $p\in \mathcal{P}_{2d}(\R^{2\times 2})$ be an isotropic polynomial. Then there are coefficients $C=(c_{ij})_{0\le i+j\le d}\in \R^{\binom{d+2}{2}}$ such that
    \begin{equation}\label{eq: isotropic polynomials}
        p(X)=\sum_{0\le i+j\le d}c_{ij}\vert X\vert^{2i}\det(X)^j. 
    \end{equation}
\end{lemma}
\begin{proof}
    By the isotropic reduction theorem \cite[Prop.~2.1 (i)]{dacorogna2005convex}, there is a polynomial $\tilde p\in \mathcal{P}_{2d}(\R^{2})$ that is invariant under the transformations in
    $$\Pi_e(2)=\{(x,y),(y,x),(-x,-y),(-y,-x)\}$$
    such that $\tilde{p}(x,y)=p(\diag(x,y))$ for all $x,y\in \R$. 
    In particular, the invariance under interchanging diagonal entries forces $\tilde{p}$ to be symmetric. By the characterization of symmetric polynomials \cite[Eq.~2.4]{macdonald1998symmetric}, there exists a polynomial $\hat{p}\in \mathcal{P}_{2d}(\R^{2})$ such that 
    $$\tilde{p}(x,y)=\hat{p}(x+y,xy).$$
    By the invariance under reflection in $\Pi_e(2)$, $\hat{p}$ must be even in the first argument. Thus, there is a polynomial $p'\in \mathcal{P}_{2d}(\R^{2})$ such that 
    $$\tilde{p}(x,y)=p'((x+y)^2,xy). $$
    The terms $2xy$ from $(x+y)^2$ can be absorbed into the second argument of $p'$, yielding the representation of $p$ by $\vert X\vert^2$ and $\det X$ in \eqref{eq: isotropic polynomials}.
\end{proof}

Although $p$ is generated by simple terms $\vert X\vert^2$ and $\det X$, it does not follow that each summand in the SOS representation of the Hessian is of this form. For instance, consider the isotropic (rank-one) convex function $p(X)=\tfrac{1}{2}\vert X\vert^2$. Then $D^2p(A)[H,H]=\vert H\vert^2$. A sum-of-squares representation is $\sum_{i=1}^4H_i^2$, however, none of these summands is isotropic. An invariant polynomial can therefore require noninvariant summands in its SOS representation

Instead of the singular value reduction, we can take a closer look at the SOS matrix $\G$ and study the additional symmetries that isotropy forces upon it. As it turns out, the matrix can be chosen to have a block-diagonal structure. 

Let us define a group action for $R,S\in \mathrm{SO}(2)$, $A,H\in \R^{2\times 2}$ via
\begin{equation}
    (R,S)\cdot (A,H):= (RAS^T,RHS^T).\label{eq: Simulataneous Rotations}
\end{equation}
\begin{lemma}\label{lem: block diagonal decomposition}
    Let $d\in \N_{>0}$ be an integer. Then, there exist an index set $\Omega_+$, natural numbers $m_0\in \N$,$m_\omega\in \N_{>0}$, $\omega\in \Omega_+$ with $m_0+\sum_{\omega\in \Omega_+}2m_\omega=N_d$, and an invertible matrix $\mathbb{F}\in \R^{N_d\times N_d}$ that only depend on $d$, such that for every for every polynomial $p$ invariant under the action \eqref{eq: Simulataneous Rotations} and admitting a representation $p(A,H)=z_d(A,H)^T\G z_d(A,H)$ with $\G\succeq0$, there exist positive semidefinite matrices
    $\mathbb K_\omega \in \R_{\mathrm{Sym}}^{2m_\omega\times 2m_\omega}$, $\mathbb K_0\in \R_{\mathrm{Sym}}^{m_0\times m_0}$ such that  $$p(A,H)=z_d(A,H)^T\mathbb{F}^T\biggl(\mathbb{K}_0\oplus \bigoplus_{\omega\in \Omega_+}\mathbb{K}_\omega\biggr)\mathbb{F}z_d(A,H).$$
    Moreover, $\mathbb{F}$ has an explicit form, given at the end of the proof below.
\end{lemma}
\begin{proof}
    The proof relies on a complexification argument, let us therefore write $z_d^{\mathbb R}(A,H)
=(A^\alpha H_{ij})_
{\substack{|\alpha|\le d-1\\i,j=1,2}}$ to make clear that this concerns real monomials. For $2\times 2$-matrices, instead of applying the matrix to a real vector, we can apply it to the complex numbers $z,\bar{z}$ via $A(z)=az+b\bar z.$
We will write
\begin{align*}
a&=\frac12(A_{11}+A_{22}+i(A_{21}-A_{12})),\\
b&=\frac12(A_{11}-A_{22}+i(A_{21}+A_{12})),\\
\xi &=\frac12(H_{11}+H_{22}+i(H_{21}-H_{12})),\\
\eta &=\frac12(H_{11}-H_{22}+i(H_{21}+H_{12})),
\end{align*}
and define a complex monomial vector
$$z_d^{\C}(A,H) =(a^r\bar a^s b^t\bar b^u\zeta)_{\substack{r+s+t+u\le d-1\\ \zeta \in\{\xi,\bar\xi,\eta,\bar\eta\}}}.$$
We denote the complex linear transform $A\to (a,\bar{a},b,\bar{b})$ by $C$, and the respective transformation on the full monomial vector by $C_d$. Then, for a positive semidefinite $\G$, 
\begin{align*}
    p(A,H)&=z_d^{\R}(A,H)^*\G^{\R}z_d^{\R}(A,H)=z_d^{\C}(A,H)^*(C_d^{-1})^* \G^{\R}C_d^{-1}z_d^{\C}(A,H).
\end{align*}Let us set $\G^{\C}=(C_d^{-1})^* \G^{\R}C_d^{-1}$, which is a positive semidefinite Hermitian matrix. 
Now, let there be given two rotations $R_\phi,R_\theta$ with respective angles. Then
\begin{align}\label{eq: Phase of bi-rotation}
    C(R_\phi AR_\theta^T)=(e^{i(\phi-\theta)}a,e^{i(-\phi+\theta)}\bar{a}, e^{i(\phi+\theta)}b,e^{-i(\phi+\theta)}\bar{b}).
\end{align}
For a pair of angles $(\phi,\theta)$, we write $\omega_w=(\pm 1,\pm 1)$ for the signs of the angles, i.e.,
$$\omega_{a}= (1,-1),\ \omega_{\bar{a}}= (-1,1),\ \omega_{b}= (1,1),\ \omega_{\bar{b}}= (-1,-1),$$
and analogously for $\xi,\eta$. 
For a complex monomial, we define its label
$$\omega(a^r\bar{a}^sb^t\bar{b}^{u}\zeta):=(r-s+t-u,-r+s+t-u)+\omega_\zeta\in \mathbb{Z}^2. $$
Let $z_\omega$ be the subvector of $z^{\C}_d$ of all monomials with label $\omega$ and denote the length of this vector by $m_\omega$. In particular, $\sum_{\omega\in \mathbb{Z}^2}m_\omega=N_d$. Let $z_0$ correspond to the zero label. Set $\Omega:=\{\omega\in\mathbb Z^2\mid m_\omega>0\}$. Conjugation pairs each monomial $a^r\bar a^s b^t\bar b^u\zeta$ of weight $\omega$ with $a^s\bar a^r b^u\bar b^t\bar\zeta$ of weight $-\omega$.
Choose $\Omega_+$ to contain exactly one representative of each nonzero pair $\{\omega,-\omega\}$. Fix an order on $\Omega_+$. We can then find a permutation $\mathbb{P}$ such that 
$$\tilde{z}_d:=\mathbb{P}z_d^{\C}=(z_\omega)_{\omega\in \Omega}, \quad \text{and}\quad z_{-\omega}=\overline{z_\omega}.$$
Then set $\tilde{\G}:=(\mathbb{P}^{-1})^*\G^{\C}\mathbb{P}^{-1}$, i.e., $p={\tilde{z}_d}^*\tilde{\G}\tilde{z}_d$. 
Now, let us compute $z_\omega$ for bi-rotated monomials. We get
\begin{align*}
    z_\omega(R_\phi AR_\theta^T,R_\phi HR_\theta^T)=e^{i(\omega_1\phi+\omega_2\theta)}z_\omega(A,H),
\end{align*}hence the action operates diagonally on each $z_\omega$. Let us write $T(\phi,\theta):=(e^{i(\omega_1\phi+\omega_2\theta)}\mathbb{I}_{m_\omega})_{\omega\in\Omega}$. By the invariance of $p$, we can replace $\tilde{\G}$ by 
$$\hat{\G}:=\frac{1}{4\pi^2}\int_{[0,2\pi)\times [0,2\pi)}T(\phi,\theta)^* \tilde{\G}T(\phi,\theta)\dd (\phi,\theta).$$
Take a block $\hat{\G}_{\omega,\chi}$ corresponding to different labels. Then 
\begin{align*}
    \hat{\G}_{\omega,\chi}:=\frac{1}{4\pi^2}\int_0^{2\pi}\int_0^{2\pi} e^{i((\chi_1-\omega_1)\phi+
(\chi_2-\omega_2)\theta)}\dd \phi\dd \theta
\tilde \G_{\omega,\chi}=\begin{cases}
    \tilde \G_{\omega,\chi} & \text{if } \omega=\chi,\\
    0 & \text{otherwise}.
\end{cases}
\end{align*}
In particular, $\hat{\G}$ is a block-diagonal matrix with Hermitian positive semidefinite blocks $\G_{\omega,\omega}$. By further averaging the block matrix with its conjugate-reflected counterpart, we can ensure that $\G_{-\omega,-\omega}=\overline{\G_{\omega,\omega}}$. Since $\hat{\G}$ is Hermitian, it follows that 
$$z_{-\omega}^*\hat{\G}_{-\omega,-\omega}z_{-\omega}=\overline{z_\omega}^* \overline{\hat{\G}_{\omega,\omega}}\overline{z_\omega}=z_\omega^*\hat{\G}_{\omega,\omega} z_\omega.$$
It remains to transfer the decomposition back to the original real vector space. Write
$$A_\omega=\operatorname{Re}\hat{\G}_{\omega,\omega}, \qquad B_\omega=\operatorname{Im}\hat{\G}_{\omega,\omega}.$$
Since $\hat{\G}$ is Hermitian, $A_\omega$ is symmetric and $B_\omega$ is skew-symmetric. The realification is given by
$$\mathbb{K}_\omega:=\mathcal R(\hat{\G}_{\omega,\omega})= \begin{pmatrix}
A_\omega & -B_\omega\\ 
B_\omega & A_\omega
\end{pmatrix}.$$
Setting $w_\omega=\sqrt{2}\begin{pmatrix} \operatorname{Re}z_\omega\\ \operatorname{Im}z_\omega \end{pmatrix}=:L_\omega(z_\omega)$, we obtain
$$z_\omega^*\hat{\G}_\omega z_\omega+ z_{-\omega}^*\hat{\G}_{-\omega}z_{-\omega} =w_\omega^T\mathbb K_\omega w_\omega.$$
For the zero weight, no individual monomial $q$ is equal to its conjugate $\bar{q}$, since conjugation always exchanges $\zeta$ for $\bar{\zeta}$. We can thus extract the real and imaginary parts as $\operatorname{Re}(q)=\frac{q+\bar{q}}{2}$, and $\operatorname{Im}(q)=\frac{q-\bar{q}}{2i}$
$$w_0=(\operatorname{Re}q_1,\operatorname{Im}q_1,\ldots,
\operatorname{Re}q_s,\operatorname{Im}q_s)^T.$$
This gives a linear isomorphism $L_0$ from the complex monomials in $z_0$ to real monomials given by their real and imaginary parts $w_0$. Conjugating $\G_{00}$ with $L_0^{-1}$, we obtain a real $\mathbb{K}_0$. 
We therefore have a decomposition 
$$p(A,H)=w_d(A,H)^T\biggl(\mathbb{K}_0\oplus \bigoplus_{\omega\in \Omega_+}\mathbb{K}_\omega\biggr)w_d(A,H).$$
Lastly, the linear isomorphism $\mathbb{F}$ between $z^{\R}_d$ to $w_d$ is given as constructed: first, we map $z^{\R}_d$ to a vector of complex monomials using $C_d$, second, we permute using $\mathbb{P}$, and third, we realify using $L_0$ and $L_\omega$. 
\end{proof}
Note that by construction
$$\mathbb K_\omega J_\omega = J_\omega\mathbb K_\omega, \qquad J_\omega=
\begin{pmatrix} 
0&-I\\
I&0
\end{pmatrix}. $$

Finally, we can combine the decomposition of the SOS matrix and the generation of $p$ by the two elementary functions and obtain a reduced semidefinite program for the isotropic Morrey problem. 
Let
$$\mathcal A_d := \{\alpha\in\mathbb N_0^4:1\leq|\alpha|\leq 2d-2\}, \qquad \mathcal I_d :=\{(i,j)\in\mathbb N_0^2:1\leq i+j\leq d-1\}.$$
The Hessian decomposition can be averaged over the rotation group, thus we restrict the determinant multiplier to isotropic polynomials of the form
$$\Lambda_\mu(X) = \sum_{(i,j)\in\mathcal I_d}
\mu_{ij}|X|^{2i}\det(X)^j.$$
Let $E_d\in\mathbb R^{\mathcal A_d\times\mathcal I_d}$ be the coefficient matrix of this basis.
Thus, if $\lambda$ denotes the usual monomial coefficient vector of $\Lambda_\mu$, then $\lambda=E_d\mu$, and we further set $\beta_u:=E_d^Tb_u$. 
\begin{corollary}
    For any $u\in C^\infty(\mathbb{T}^2,\R^2)$ and $d\in \N$, we let $\mathbb{M}_u\in \R^{N_d\times N_d}$ and $b_u\in \R^{\binom{2d+2}{4}-1}$, be defined as in Thm.~\ref{Thm: Morrey Test}. Let $\mathbb{F}$, $\Omega_+$, $m_0\in \N$,$m_\omega\in \N_{>0}$, be as in Lem.~\ref{lem: block diagonal decomposition}. Let  $\mathbb K_\omega \in \R_{\mathrm{Sym}}^{2m_\omega\times 2m_\omega}$, $\mathbb K_0\in \R_{\mathrm{Sym}}^{m_0\times m_0}$.
    Define $\widehat{\mathbb{M}}_u:=\mathbb{F}\mathbb{M}_u\mathbb{F}^T$.
The reduced semidefinite program is given by
\begin{align}\label{eq: SDP iso}
\begin{aligned}
\rho_d^{\mathrm{iso}}(u) =& \inf_{\mathbb K,\mu}
\biggl[ \langle\widehat{\mathbb M}_{u,00},\mathbb K_0\rangle + \sum_{\omega\in\Omega_+} \langle\widehat{\mathbb M}_{u,\omega,\omega}, \mathbb K_\omega\rangle \biggr]\\
\text{subject to }  &\mathscr C_d\big(
\mathbb F^T\mathbb K\mathbb F,E_d\mu
\big)=0, \ \beta_u^T\mu=-1,\\
&\mathbb K_0\succeq0,\ \mathbb K_\omega\succeq0, \ \mathbb K_\omega \ J_\omega
=J_\omega\mathbb K_\omega.
\end{aligned} \tag{SDP-iso}
\end{align}
    Then there is an isotropic rank-one convex function on $\R^{2\times 2}$ failing to be quasiconvex at $0$ with witness $u$ if and only if there is a $d\in \N$ such that $\rho^{\mathrm{iso}}_d(u)<1$.
\end{corollary}
\begin{Bem}
The same criterion holds at any prescribed $A\in\mathbb R^{2\times2}$. Define
\begin{align*}
\mathbb M_{u,A} &=\int_{\mathbb T^2}\int_0^1(1-t) z_d(A+tDu,Du)\otimes z_d(A+tDu,Du)\dd t\dd x,\\
(b_{u,A})_\alpha &=\int_{\mathbb T^2}\int_0^1(1-t)
(A+tDu)^\alpha\det(Du)\dd t\dd x .
\end{align*}
Taylor’s formula at $A$ and $\int_{\mathbb T^2}Du\dd x=0$ give
$$\int_{\mathbb T^2}p(A+Du)\dd x-p(A) =\langle\mathbb M_{u,A},\G\rangle + b_{u,A}^T\lambda.$$
Hence \eqref{eq: SDP iso}, with $\widehat{\mathbb M}_{u,A}=\mathbb F\mathbb M_{u,A}\mathbb F^T$ and $\beta_{u,A}=E_d^Tb_{u,A}$, detects a violation at $A$ precisely when its value is below $1$ for some $d$. For the density implication, approximate the isotropic integrand on the compact set $\{A\}\cup(A+Du(\mathbb T^2))$.
\end{Bem}
We can further simplify this by taking $C$ as in \eqref{eq: isotropic polynomials}. Consequently, the number of compatibility constraints decreases, since there are fewer monomials to consider. 
Note that the general matrix $\G$ has $N_d^2=\mathcal{O}(d^{4\cdot 2})$ entries, and its Cholesky decomposition requires $N_d^3=\mathcal{O}(d^{12})$ computations. Every $\mathbb{K}$ has only $\mathcal{O}(d^{2\cdot 2})$ entries, and there are $\mathcal{O}(d^2)$ such blocks\footnote{Since the sum of the exponents is at most $d-1$, both components of $\omega$ lie between approximately $-d$ and $d$. For fixed $\omega,\zeta$, the weight determines $r-s$ and $t-u$, leaving $\mathcal{O}(d)$ choices for each of $\min(r,s)$ and $\min(t,u)$, hence $m_\omega=\mathcal{O}(d^2)$.}. The Cholesky decomposition therefore takes $\mathcal{O}(d^{2\cdot 3}\cdot d^2)=\mathcal{O}(d^8)$ computations.

\subsection{Corollaries on the symmetric Morrey problem}

Symmetric rank-one convexity is a variant of rank-one convexity that is tailored to integrands arising in geometrically linear elasticity. Crucially, energies there depend on the symmetrized gradient $\varepsilon(u)=\tfrac{1}{2}(Du+Du^T)$. Let $f\colon \Symtwo\to \R$ be continuous. Define
$$\mathbb{E}(u):=\int_{\Omega}f(\varepsilon(u))\dd x$$
as in \eqref{eq: general energy}. We let $\pi_S\colon \R^{2\times 2}\to \Symtwo$, $A\mapsto \tfrac{1}{2}(A+A^T)$ be the natural projection, we can write this as
$$\mathbb{E}(u)=\int_\Omega f\circ \pi_S(Du)\dd x.$$
Importantly, this is now of the same form as \eqref{eq: general energy}. In particular, the relevant property for the lower semicontinuity of $\mathbb{E}$ is the quasiconvexity of $f\circ \pi_S$, which is equivalent to $f$ being symmetric quasiconvex. Analogously, one defines $f$ to be symmetric rank-one convex if $f\circ \pi_S$ is rank-one convex (see \cite{Kreisbeck19, DolGam2026,GamKre2026}). The symmetric rank-one cone consists of symmetrized rank-one directions, which can have ordinary rank two. To be more precise, a matrix $A\in \Symtwo$ is symmetric rank-one if there are $a,b\in \R^2$ such that $A=a\odot b=\frac{1}{2}(a\otimes b+b\otimes a)$. In two dimensions, this is the same as $\det A\le 0$. 

The strong symmetric Morrey problem, i.e., whether symmetric rank-one convexity implies symmetric quasiconvexity in two dimensions\footnote{In dimensions $\ge 3$, this is known to be false in general, see \cite{gammer2026symmetricmorreyproblemdimensions}}, is interesting in its own right, both for mathematical and materials science reasons, as it relates to the laminate structures of multiphase materials. However, such a counterexample $f$ in the symmetric case already yields a counterexample to the non-symmetric case by lifting it to $f\circ \pi_S$. Therefore, establishing the symmetric analog of our results can be valuable for the non-symmetric SDP. 

Let us first define symmetric SOSRC polynomials
\begin{Definition}
     A polynomial $p\in \mathcal{P}(\Symtwo)$ is \emph{symmetric sum-of-squares rank-one convex} if there are polynomials $r_1,...,r_N\in \mathcal{P}(\Symtwo \times \R^2\times\R^2)$, $N\in \N$ such that for all $A\in \Symtwo$, $a\in \R^2$, $b\in \R^2$ there holds
    \begin{equation}
         D^2p(A)[a\odot b,a\odot b]=\sum_{\nu=1}^Nr_{\nu}(A,a,b)^2.
    \end{equation}
\end{Definition}

Next, we show that these polynomials are dense in the symmetric rank-one convex functions. As the proof of Lemma \ref{lem: polynomial selection} shows, it suffices to find a polynomial that is symmetric rank-one convex on a compact set. 

\begin{Proposition}
    Let $f\in C^0(\Symtwo)$ be a symmetric rank-one convex function. Then for any compact $K\subset  \Symtwo $, and any $\delta>0$, there is a \emph{symmetric $\mathrm{SOSRC}$} polynomial $p\in   \mathcal{P}(\Symtwo)$ such that 
    $$\sup_{A\in K}\vert f(A)-p(A)\vert <\delta.$$
\end{Proposition}
\begin{proof}
    We only sketch the proof. Choose a compact set $K$, and a real number $R>0$ such that $K\Subset Q_{R/2}$. Then approximate $f$ by $C^2$-functions in the local $C^0$-topology by mollification. Choose a smooth mollification $f_\delta$ sufficiently close to $f$ and add $\delta\vert \cdot \vert^2$ to make it strictly rank-one convex. Approximate $f_\delta+\delta\vert \cdot \vert^2$ with a polynomial $p$ on $Q_R$ in $C^2$ on a scale $\delta/2$. Necessarily, $p$ is rank-one convex on $Q_R$. 

In order to apply the polynomial selection lemma \ref{lem: polynomial selection}, we cannot use Finsler's lemma in the symmetric case. Instead, we can use the characterization of the symmetric rank one cone as $\{\det \le 0\}$. Then the Yakubovich S-lemma (\cite{yakubovich1971s}, \cite[Thm.~2.2]{polik2007survey}) gives a nonnegative multiplier $\mu_A$ for every point $A$. As in lemma \ref{lem: polynomial selection}, using a partition of unity, we can then select a continuous representation of the multipliers and approximate its square root by a polynomial $q$ on $Q_R$. This yields a $\Lambda$ which is itself the square of the polynomial $q$. 

Next, we apply the matrix Lasserre--Netzer construction from lemmas \ref{lem: matrix Lasserre} and \ref{lem: Hessian SOS}, where we reduce the number of variables from $N=4$ to $N=3$ by symmetry. The rest of the proof proceeds analogously to the proof of Thm.~\ref{Thm: Density of SOSRC}, where the resulting polynomial is symmetric SOSRC globally, and approximates $f$ on $Q_{R/2}$ on a scale $\delta$. 
\end{proof}

We identify $\Symtwo\cong \R^3$ and write $\langle A,H\rangle_F=a_1h_1+a_3h_3+2a_2h_2$. The following corollary shows that the symmetric SOSRC condition forces the determinant multiplier to be a negative sum of squares.

\begin{corollary}
    Let $p\in \mathcal{P}_{2d}(\R^3)$ be symmetric \emph{SOSRC}. Then there exists an integer $N\in \N$, polynomials $q_1,...,q_N\in \mathcal{P}_{d-1}$, $v_\nu\in \mathcal{P}_{d-1}(\R^3)^3$, such that for every $A,H\in \R^3$ there hold
    \begin{equation}
        D^2p(A)[H,H]+\biggl(\sum_{\nu=1}^N q_\nu(A)^2\biggr)\det H=\sum_{\nu=1}^N \langle v_\nu(A),H\rangle_F^2.
    \end{equation}
    In particular, the multiplier $\Lambda:=\sum_{\nu=1}^N q_\nu^2$ is \emph{SOS}. 
\end{corollary}
\begin{proof}
The expression $D^2p(A)[a\odot b,a\odot b]$ is homogeneous of degree two in $a$ and separately homogeneous of degree two in $b$. Write $\sum_{\nu=1}^Nr_{\nu}(A,a,b)^2$ for its SOSRC representation.
For a sum of squares, the highest-degree homogeneous parts cannot cancel. Thus each summand $r_\nu$ is homogeneous of degree one in $a$ and $b$. In particular, we can represent $r_\nu$ via matrix polynomials $B_\nu$ 
$$r_\nu(A,a,b)=a^TB_\nu(A)b.$$
Since $D^2p(A)$ is symmetric in $a,b$, we can average the representation. Let $v_\nu$ denote the symmetric part of $B_\nu$, and $K_\nu$ denote its skew-symmetric part. Note that first 
$$K_\nu(A)=\kappa_\nu(A)\begin{pmatrix}0&1\\-1 &0\end{pmatrix}, \quad \text{where }\kappa_\nu \in \mathcal{P}_{d-1}(\R^3),$$ 
and second that $a^TK_\nu(A)b=\kappa_\nu(A)(a_1b_2-a_2b_1)$.
Then
\begin{align*}
    D^2p(A)[a\odot b,a\odot b]&=\sum_{\nu=1}^N \langle v_\nu(A),a\odot b\rangle^2_F+\sum_{\nu=1}^N \kappa_\nu(A)^2(a_1b_2-a_2b_1)^2\\
    &=\sum_{\nu=1}^N \langle v_\nu(A),a\odot b\rangle^2_F-4\sum_{\nu=1}^N \kappa_\nu(A)^2\det(a\odot b).
\end{align*}
Since this holds for all $H=a\odot b$, it holds on the set $\{\det\le 0\}$, in particular on the open set $\{\det <0\}$. Thus, the polynomial identity must hold on $\Symtwo$. Finally, set $q_\nu=2\kappa_\nu$.
\end{proof}

We get an immediate analogue to Cor.~\ref{Cor: Generation of SOSRC functions}
\begin{corollary}
Let $d\in \N$, let $m_{d-1}(X)=(X^\alpha)_{\vert \alpha\vert\le d-1}$ be monomials up to degree $d-1$ in $3$ variables. Set $z_d^{\mathrm{Sym}}(X,H)=m_{d-1}(X)\otimes H$, $n_d=\binom{d+2}{3}$. Let $\G\in \R_{\mathrm{Sym}}^{3n_d\times 3n_d}$, $\mathbb{L}\in \R_{\mathrm{Sym}}^{n_d\times n_d}$ be positive semidefinite, and let $C\in \R^{\binom{2d+3}{3}-4}$. For coefficients $C$, we set $p_C:=\sum_{\vert \alpha\vert =2}^{2d}C_\alpha X^\alpha$, and define a feasibility condition
\begin{align*}
    \mathscr{F}_d^{\mathrm{Sym}}(C,\G,\mathbb L)[X,H]&:=D^2p_C(X)[H,H]-z_d^{\mathrm{Sym}}(X,H)^T\G z_d^{\mathrm{Sym}}(X,H)\\
    & \qquad +m_{d-1}(X)^T\mathbb{L}m_{d-1}(X)\det(H)=0.
\end{align*}
Then, any polynomial $p$ on $\Symtwo$ of degree at most $2d$ is symmetric $\mathrm{SOSRC}$ if and only if it is of the form $p_C+p_{\mathrm{aff}}$ for some affine $p_{\mathrm{aff}}$, $C,\G,\mathbb L$ satisfying
$$\mathscr{F}_d^{\mathrm{Sym}}(C,\G,\mathbb L)\equiv 0, \quad \text{and} \quad \G\succeq 0.$$
\end{corollary}

Lastly, the analog to \eqref{eq: semidefinite program} reads as follows, where $\mathbb{D}_{\mathrm{Sym}}$ represents the bilinear determinant mapping in a symmetric setting, i.e., if $H=(h_{11},h_{12},h_{22})$, then $H^T\mathbb{D}_{\mathrm{Sym}} H=h_{11}h_{22}-h_{12}^2$.

\begin{corollary}
 For any $u\in C^\infty(\mathbb{T}^2,\R^2)$ and $d\in \N$, we define matrices $\mathbb{M}_u^{\mathrm{Sym}}\in \R^{3n_d\times 3n_d}$, $\mathbb{N}_u^{\mathrm{Sym}}\in \R^{n_d\times n_d}$ by
 \begin{align*}
     \mathbb M_u^{\mathrm{Sym}} &= \int_{\mathbb T^2}\int_0^1
(1-t)\,
z_d^{\mathrm{Sym}}(t\varepsilon(u),\varepsilon(u))
z_d^{\mathrm{Sym}}(t\varepsilon(u),\varepsilon(u))^T
\dd t\dd x\\
 \mathbb N_u^{\mathrm{Sym}} &= \int_{\mathbb T^2}\int_0^1 (1-t)\det \varepsilon(u) m_{d-1}(t\varepsilon(u))m_{d-1}(t\varepsilon(u))^T \dd t\dd x.
 \end{align*}
   Let $\mathcal{C}_d^{\mathrm{Sym}},\mathcal{H}_d^{\mathrm{Sym}}$ be the symmetric analogs of the operators in \eqref{eq: Linear Operators} and define
$$\mathcal{L}_d^{\mathrm{Sym}}\colon \R^{{n_d}\times {n_d}}_{\mathrm{Sym}}\to \Symm, \ \mathbb{L}\mapsto m_{d-1}(X)^T\mathbb{L}m_{d-1}(X)\mathbb{D}_{\mathrm{Sym}}.$$
   We define a semidefinite program via
    \begin{align}\label{eq: symmetric semidefinite program}\tag{Sym SDP}
        \begin{aligned}
&\rho_d^{\mathrm{Sym}}(u)=\inf_{\G,\mathbb{L},C}\, \langle \mathbb{M}_u^{\mathrm{Sym}},\G\rangle \\
            &\mathrm{subject\ to} \ \mathcal{H}^{\mathrm{Sym}}(p_C)=\mathcal{C}_d^{\mathrm{Sym}}(\G)-\mathcal{L}_d^{\mathrm{Sym}}(\mathbb{L}), \ \langle\mathbb N_u^{\mathrm{Sym}},\mathbb{L}\rangle= 1, \ \mathrm{and}\ \G,\mathbb{L}\succeq 0.
        \end{aligned}
    \end{align}
    Then there is a symmetric rank-one convex function on $\Symtwo$ failing to be symmetric quasiconvex at $0$ with witness $u$ if and only if there is a $d\in \N$ such that $\rho_d^{\mathrm{Sym}}(u)<1$.
\end{corollary}
\begin{proof}
We calculate
    \begin{align*}
 \int_{\mathbb T^2}p_C(\varepsilon(u))\dd x&=  \int_{\mathbb T^2}\int_0^1 (1-t)D^2p_C(t\varepsilon(u))[\varepsilon(u),\varepsilon(u)]\dd t\dd x\\
 &= \langle\mathbb M_u^{\mathrm{Sym}},\G\rangle - 
\langle\mathbb N_u^{\mathrm{Sym}},\mathbb{L}\rangle.
\end{align*}
The corollary follows by the same argument as in Thm.~\ref{Thm: Morrey Test}.
\end{proof}
  The symmetric SDP consists of two positive semidefinite matrix variables of orders $3n_d$ and $n_d$ respectively, where $n_d=\binom{d+2}{3}=\mathcal{O}(d^3)$. Their Cholesky factorizations thus require $\mathcal{O}((3n_d)^3+n_d^3)=\mathcal{O}(d^9)$ operations, compared to $\mathcal{O}(d^{12})$ for the unrestricted SDP.

\section{Outlook}

Although we were unable to construct a rank-one convex function that is not quasiconvex in $\R^{2\times 2}$, Thm.~\ref{Thm: Density of SOSRC} and our SDP formulation show that for any prescribed test function $u$ causing a violation of the quasiconvexity inequality for some $f$, a sufficiently high polynomial degree will detect the violation within the SOSRC class. 

Thus, a potential route for future work could be to improve and speed up the SDP algorithm. This can potentially be done in several ways. A first one is to limit the computations to subclasses with invariances as in Section \ref{Sec: Corollaries}, or even additional symmetries. A second route could be to incorporate dimension reductions, parallelize searches across test functions $u$, as well as extend the program to higher degrees $d$ using more computing resources. 
As a third option, since these numerical schemes still involve an $\mathcal{O}(d^{r})$ operations in general (where $r\ge 8$), future work could also focus on exchanging the optimization algorithm for a Monte-Carlo procedure or variants of gradient descent. Even if these do not obtain the optimal value in the SDP, they might still be able to detect a polynomial where the objective is less than $1$. 

On the theoretical side, if rank-one convexity indeed implies quasiconvexity, then Thm.~\ref{Thm: Density of SOSRC} shows that it suffices to establish this result for SOSRC polynomials.
The representation theorem $\ref{Thm: Representation}$ gives rise to a novel approach to such a proof by showing that for any $\Lambda$, every compatible $\G$ leads to a quasiconvex $p_C$. This could, for example, be done by proving a lower bound on the SDP objective (or equivalently a dual feasible value $\eta\ge 1$).

To disprove the implication, it would suffice to construct suitable $\Lambda$ and $u$, and then to show that there is a degree $d$ at which the admissible matrices $\G$ admit a quasiconvexity violation witness.

\subsubsection*{Acknowledgments}
The research for this paper has been conducted at the University of Regensburg and was brought to an end during a research stay of the author at the University of Oxford. The author thanks Jan Kristensen and Gabriele Cassese for many insightful conversations.

This work was funded by the Deutsche Forschungsgemeinschaft (DFG). We gratefully acknowledge the support within the project SPP 2256
(``Variational Methods for Predicting Complex Phenomena in Engineering Structures
and Materials''), project ID 441468770.

\subsubsection*{AI disclosure}

The original draft of section 2 only concerned the $2\times 2$-case. After review of this draft with ChtGPT 6, the AI suggested the extension to the $\min\{m,n\}=2$ case using Calderón's theorem instead of Finsler's lemma. Also in this original draft, the proof of the density results used the polynomial density argument by Heinz and then approximated the rank-one convex polynomials with SOSRC polynomials. The AI pointed out that this was unnecessary, since for the application of Lemma~\ref{lem: polynomial selection}, it suffices to have a polynomial to be rank-one convex on $Q_R$. Such polynomial approximates can easily be obtained by the Weierstrass theorem. In this sense, the arguments here prove a stronger version than Heinz' without relying on the latter in the first place. 

Besides this, the proofs in section 2, 3, 4, and 5.2 have been conducted without AI use. In section 5.1, the density result was obtained by the author. The proof of the block-diagonal structure was originally proposed by the AI, the author then made this rigorous and (hopefully) human-readable. 

Other than that, Claude Sonnet 5.5 was used for the first implementation of the numerical scheme. The author then improved the run-time and stability properties of the code. As described in the numerics section, the run time of the code still provides a challenge for basically any hardware and only allowed us to perform experiments up to degree $2d=12$ so far.

The AI models are not authors of this work. This paper has been written and checked entirely by the author, who takes full responsibility for its contents.

\printbibliography[heading=bibintoc, title={References}]

@article{morrey1952quasi,
  author = {Morrey, Jr., Charles B.},
  title = {Quasi-convexity and the lower semicontinuity of multiple integrals},
  journal = {Pacific Journal of Mathematics},
  volume = {2},
  number = {1},
  pages = {25--53},
  year = {1952}
}

@article{vsverak1992rank,
  author = {{\v{S}}ver{\'a}k, Vladim{\'\i}r},
  title = {Rank-one convexity does not imply quasiconvexity},
  journal = {Proceedings of the Royal Society of Edinburgh Section A: Mathematics},
  volume = {120},
  number = {1--2},
  pages = {185--189},
  year = {1992}
}

@article{kristensen1999non,
  author = {Kristensen, Jan},
  title = {On the non-locality of quasiconvexity},
  journal = {Annales de l'Institut Henri Poincar{\'e} C, Analyse non lin{\'e}aire},
  volume = {16},
  number = {1},
  pages = {1--13},
  year = {1999}
}

@article{guerra2021numerical,
  author = {Guerra, Andr{\'e} and Teixeira da Costa, Rita},
  title = {Numerical evidence towards a positive answer to {Morrey}'s problem},
  journal = {Revista Matem{\'a}tica Iberoamericana},
  volume = {38},
  number = {2},
  pages = {601--614},
  year = {2022}
}

@article{grabovsky2019explicit,
  author = {Grabovsky, Yury and Truskinovsky, Lev},
  title = {Explicit relaxation of a two-well {Hadamard} energy},
  journal = {Journal of Elasticity},
  volume = {135},
  number = {1--2},
  pages = {351--373},
  year = {2019}
}

@article{grabovsky2019rank,
  author = {Grabovsky, Yury and Truskinovsky, Lev},
  title = {When rank-one convexity meets polyconvexity: an algebraic approach to elastic binodal},
  journal = {Journal of Nonlinear Science},
  volume = {29},
  number = {1},
  pages = {229--253},
  year = {2019}
}

@article{chenchiah2008relaxation,
  author = {Chenchiah, Isaac V. and Bhattacharya, Kaushik},
  title = {The relaxation of two-well energies with possibly unequal moduli},
  journal = {Archive for Rational Mechanics and Analysis},
  volume = {187},
  number = {3},
  pages = {409--479},
  year = {2008}
}

@article{KristensenGuerra2022,
  author = {Guerra, Andr{\'e} and Kristensen, Jan},
  title = {Automatic quasiconvexity of homogeneous isotropic rank-one convex integrands},
  journal = {Archive for Rational Mechanics and Analysis},
  volume = {245},
  pages = {479--500},
  year = {2022}
}

@article{Kreisbeck19,
  author = {Omar Boussaid and Carolin Kreisbeck and Anja Schl{\"o}merkemper},
  title = {Characterizations of symmetric polyconvexity},
  journal = {Archive for Rational Mechanics and Analysis},
  volume = {234},
  pages = {417--451},
  year = {2019}
}

@article{martin2017rank,
  author = {Martin, Robert J. and Ghiba, Ionel-Dumitrel and Neff, Patrizio},
  title = {Rank-one convexity implies polyconvexity for isotropic, objective and isochoric elastic energies in the two-dimensional case},
  journal = {Proceedings of the Royal Society of Edinburgh Section A: Mathematics},
  volume = {147},
  number = {3},
  pages = {571--597},
  year = {2017}
}

@article{grabovsky2018microstructure,
  author = {Grabovsky, Yury},
  title = {From microstructure-independent formulas for composite materials to rank-one convex, non-quasiconvex functions},
  journal = {Archive for Rational Mechanics and Analysis},
  volume = {227},
  number = {2},
  pages = {607--636},
  year = {2018}
}

@misc{cassese2026solutionmorreysproblemmathbbr2times,
  author = {Gabriele Cassese},
  title = {A solution to {Morrey}'s problem in {$\mathbb{R}^{2\times m}$}},
  year = {2026},
  eprint = {2608.03488},
  archiveprefix = {arXiv},
  primaryclass = {math.AP},
  url = {https://arxiv.org/abs/2608.03488}
}

@misc{agazzi2026morreysproblemmathbbr2times,
  author = {Andrea Agazzi and Giuseppe Bruno and Andr{\'e} Guerra and Federico Pasqualotto},
  title = {{Morrey}'s problem in {$\mathbb{R}^{2\times 4}$} and {$\mathbb{R}^{3\times 3}_{\mathrm{sym}}$}},
  year = {2026},
  eprint = {2608.12298},
  archiveprefix = {arXiv},
  primaryclass = {math.AP},
  url = {https://arxiv.org/abs/2608.12298}
}

@article{conti2003polyconvexity,
  author = {Conti, Sergio and De Lellis, Camillo and M{\"u}ller, Stefan and Romeo, Mario},
  title = {Polyconvexity equals rank-one convexity for connected isotropic sets in {$\mathbb{M}^{2\times 2}$}},
  journal = {Comptes Rendus. Math{\'e}matique},
  volume = {337},
  number = {4},
  pages = {233--238},
  year = {2003}
}

@article{harris2018two,
  author = {Harris, Terence L. J. and Kirchheim, Bernd and Lin, Chun-Chi},
  title = {Two-by-two upper triangular matrices and {Morrey}'s conjecture},
  journal = {Calculus of Variations and Partial Differential Equations},
  volume = {57},
  number = {3},
  eid = {73},
  year = {2018}
}

@article{muller1999rank,
  author = {M{\"u}ller, Stefan},
  title = {Rank-one convexity implies quasiconvexity on diagonal matrices},
  journal = {International Mathematics Research Notices},
  volume = {1999},
  number = {20},
  pages = {1087--1095},
  year = {1999}
}

@article{faraco2003quasiconvex,
  author = {Faraco, Daniel and Zhong, Xiao},
  title = {Quasiconvex functions and {Hessian} equations},
  journal = {Archive for Rational Mechanics and Analysis},
  volume = {168},
  number = {3},
  pages = {245--252},
  year = {2003}
}

@unpublished{DolGam2026,
  author = {Dolzmann, Georg and Gammer, Florian},
  title  = {On symmetric polyconvexification and its application towards a double-well problem in 2d},
  note   = {In preparation},
  year   = {2026}
}

@unpublished{GamKre2026,
  author = {Gammer, Florian and Kreisbeck, Carolin},
  title  = {Representations of symmetric polyconvex envelopes in $3$d via convex duality},
  note   = {In preparation},
  year   = {2026}
}

@misc{ConDolDon2026,
      title={A frame-indifferent rank-one convex function in three dimensions which is not quasiconvex}, 
      author={Sergio Conti and Georg Dolzmann and Patrick Dondl},
      year={2026},
      eprint={2610.31886},
      archivePrefix={arXiv},
      primaryClass={math.AP},
      url={https://arxiv.org/abs/2610.31886}, 
}

@phdthesis{Heinz2008Preservation,
author = {Heinz, Sebastian},
title = {Preservation of quasiconvexity and quasimonotonicity in polynomial approximation of variational problems},
school = {Humboldt-Universität zu Berlin, Mathematisch-Naturwissenschaftliche Fakultät II},
year = {2008},
doi = {http://dx.doi.org/10.18452/15808}
}

@article{lassere2007sos,
  title={SOS approximations of nonnegative polynomials via simple high degree perturbations},
  author={Lasserre, Jean B. and Netzer, Tim},
  journal={Mathematische Zeitschrift},
  volume={256},
  number={2},
  pages={99--112},
  year={2007},
  publisher={Springer}
}

@article{CIMPRIC2015239,
title = {Finsler's Lemma for matrix polynomials},
journal = {Linear Algebra and its Applications},
volume = {465},
pages = {239-261},
year = {2015},
author = {Jaka Cimpri{\v c}},
}

@article{blekherman2006there,
  title={There are significantly more nonegative polynomials than sums of squares},
  author={Blekherman, Grigoriy},
  journal={Israel Journal of Mathematics},
  volume={153},
  number={1},
  pages={355--380},
  year={2006},
  publisher={Springer}
}

@misc{gammer2026symmetricmorreyproblemdimensions,
      title={On the symmetric Morrey problem in three dimensions}, 
      author={Florian Gammer},
      year={2026},
      eprint={2610.35216},
      archivePrefix={arXiv},
      primaryClass={math.AP},
      url={https://arxiv.org/abs/2610.35216}, 
}

@article{GUTIERREZ20073177,
title = {An optimization algorithm applied to the Morrey conjecture in nonlinear elasticity},
journal = {International Journal of Solids and Structures},
volume = {44},
number = {10},
pages = {3177-3186},
year = {2007},
author = {Sergio Gutiérrez and Jorge Villavicencio}
}

@misc{guerra2026quasiconvexityburkholderfunctionsymmetric,
      title={Quasiconvexity of the Burkholder function on symmetric matrices}, 
      author={André Guerra},
      year={2026},
      eprint={2608.23388},
      archivePrefix={arXiv},
      primaryClass={math.AP},
      url={https://arxiv.org/abs/2608.23388}
}

@article{kirchheimautomatic,
  title={Automatic convexity of rank-1 convex functions.(English, French summaries)},
  author={Kirchheim, Bernd and Kristensen, Jan},
 journal = {C. R. Acad. Sci. Paris},
series = {I 349},
pages = {407–409.},
year = {2011}
}

@article{motzkin1967arithmetic,
  title={The arithmetic-geometric inequality},
  author={Motzkin, Theodore Samuel},
  journal={Inequalities (Proc. Sympos. Wright-Patterson Air Force Base, Ohio, 1965)},
  volume={205},
  pages={54},
  year={1967}
}

@inproceedings{artin1927zerlegung,
  title={{\"U}ber die Zerlegung definiter Funktionen in Quadrate},
  author={Artin, Emil},
  booktitle={Abhandlungen aus dem Mathematischen Seminar der Universit{\"a}t Hamburg},
  volume={5},
  number={1},
  pages={100--115},
  year={1927},
  organization={Springer}
}

@article{hilbert1888darstellung,
  title={{\"U}ber die Darstellung definiter Formen als Summe von Formenquadraten},
  author={Hilbert, David},
  journal={Mathematische Annalen},
  volume={32},
  number={3},
  pages={342--350},
  year={1888},
  publisher={Springer}
}

@article{vandenbergheisemidefinite,
  title={Semidefinite Progamming},
  author={Lieven Vanderberghe and Stephen Boyd},
  journal={SIAM Review},
  year={1996},
  month={3},
  pages={49--95},
  volume={38},
  number={1}
}

@article{ahmadi2013np,
  title={NP-hardness of deciding convexity of quartic polynomials and related problems},
  author={Ahmadi, Amir Ali and Olshevsky, Alex and Parrilo, Pablo A and Tsitsiklis, John N},
  journal={Mathematical programming},
  volume={137},
  number={1},
  pages={453--476},
  year={2013},
  publisher={Springer}
}

@article{lasserre2007sum,
  title={A sum of squares approximation of nonnegative polynomials},
  author={Lasserre, Jean B},
  journal={SIAM review},
  volume={49},
  number={4},
  pages={651--669},
  year={2007},
  publisher={SIAM}
}

@article{berg1976positive,
  title={Positive definite functions on abelian semigroups},
  author={Berg, Christian and Christensen, Jens Peter Reus and Ressel, Paul},
  journal={Mathematische Annalen},
  volume={223},
  number={3},
  pages={253--274},
  year={1976},
  publisher={Springer}
}

@article{finsler1936vorkommen,
  title={{\"U}ber das Vorkommen definiter und semidefiniter Formen in Scharen quadratischer Formen},
  author={Finsler, Paul},
  journal={Commentarii Mathematici Helvetici},
  volume={9},
  number={1},
  pages={188--192},
  year={1936},
  publisher={Springer}
}

@article{calderon1973note,
  title={A note on biquadratic forms},
  author={Calder{\'o}n, Alberto Pedro},
  journal={Linear Algebra and its Applications},
  volume={7},
  number={2},
  pages={175--177},
  year={1973},
  publisher={Elsevier}
}

@article{voss2022numerical,
  title={Numerical approaches for investigating quasiconvexity in the context of Morrey’s conjecture},
  author={Voss, Jendrik and Martin, Robert J and Sander, Oliver and Kumar, Siddhant and Kochmann, Dennis M and Neff, Patrizio},
  journal={Journal of Nonlinear Science},
  volume={32},
  number={6},
  pages={77},
  year={2022},
  publisher={Springer}
}

@article{bandeira2016role,
  title={The role of non-negative polynomials for rank-one convexity and quasi convexity},
  author={Bandeira, Luis and Pedregal, Pablo},
  journal={Journal of Elliptic and Parabolic Equations},
  volume={2},
  number={1},
  pages={27--36},
  year={2016},
  publisher={Springer}
}

@inproceedings{bandeira2009finding,
  title={Finding new families of rank-one convex polynomials},
  author={Bandeira, Lu{\'\i}s and Pedregal, Pablo},
  booktitle={Annales de l'IHP Analyse non lin{\'e}aire},
  volume={26},
  number={5},
  pages={1621--1634},
  year={2009}
}

@book{Bhatia1997MatrixAnalysis,
  author    = {Bhatia, Rajendra},
  title     = {Matrix Analysis},
  series    = {Graduate Texts in Mathematics},
  volume    = {169},
  publisher = {Springer},
  address   = {New York},
  year      = {1997}
}

@article{balakrishnan2003semidefinite,
  title={Semidefinite programming duality and linear time-invariant systems},
  author={Balakrishnan, Venkataramanan and Vandenberghe, Lieven},
  journal={IEEE Transactions on Automatic Control},
  volume={48},
  number={1},
  pages={30--41},
  year={2003},
  publisher={IEEE}
}

@incollection{ball1987does,
  title={Does rank-one convexity imply quasiconvexity?},
  author={Ball, John M},
  booktitle={Metastability and incompletely posed problems},
  pages={17--32},
  year={1987},
  publisher={Springer}
}

@article{yakubovich1971s,
  title={S-procedure in nolinear control theory},
  author={Yakubovich, Vladimir A},
  journal={Vestnik Leninggradskogo Universiteta, Ser. Matematika},
  pages={62--77},
  year={1971}
}

@article{polik2007survey,
  title={A survey of the S-lemma},
  author={P{\'o}lik, Imre and Terlaky, Tam{\'a}s},
  journal={SIAM review},
  volume={49},
  number={3},
  pages={371--418},
  year={2007},
  publisher={SIAM}
}

@book{macdonald1998symmetric,
  title={Symmetric functions and Hall polynomials},
  author={Macdonald, Ian Grant},
  year={1998},
  publisher={Oxford university press}
}

@inproceedings{dacorogna2005convex,
  title={Convex SO(N)$\times$ SO(n)-invariant functions and refinements of Von Neumann’s inequality},
  author={Dacorogna, Bernard and Mar{\'e}chal, Pierre},
  booktitle={Annales de la Facult{\'e} des Sciences de Toulouse},
  year={2005}
}

@misc{henrion2026semidefiniterelaxationsnonlinearelasticity,
      title={Semidefinite relaxations for nonlinear elasticity with energies convex in the Cauchy-Green strain tensor}, 
      author={Didier Henrion and Milan Korda and Martin Kružík and Karolína Sehnalová},
      year={2026},
      eprint={2604.13566},
      archivePrefix={arXiv},
      primaryClass={math.OC},
      url={https://arxiv.org/abs/2604.13566}
}

@misc{fantuzzi2026polyconvexitymomentssumssquares,
      title={Polyconvexity with Moments and Sums of Squares}, 
      author={Giovanni Fantuzzi and Didier Henrion and Martin Kru{ž}ík and Ajay Murali and Stephan Weis},
      year={2026},
      eprint={2604.11124},
      archivePrefix={arXiv},
      primaryClass={math.OC},
      url={https://arxiv.org/abs/2604.11124}
}

@article{mehrotra1992implementation,
  title={On the implementation of a primal-dual interior point method},
  author={Mehrotra, Sanjay},
  journal={SIAM Journal on optimization},
  volume={2},
  number={4},
  pages={575--601},
  year={1992},
  publisher={SIAM}
}

@article{choi1975positive,
  title={Positive semidefinite biquadratic forms},
  author={Choi, Man-Duen},
  journal={Linear Algebra and its applications},
  volume={12},
  number={2},
  pages={95--100},
  year={1975},
  publisher={Elsevier}
}

\end{document}